\documentclass[12pt]{article}
\usepackage{amssymb,amsmath}
\usepackage{mathrsfs}
\usepackage{graphicx}

\newtheorem{Thm}{Theorem}
\newtheorem{Lm}{Lemma}
\newtheorem{Def}{Definition}

\begin{document}

\begin{center}
\large{\bf CLASSIFICATION OF CATALYTIC BRANCHING PROCESSES AND
STRUCTURE OF THE CRITICALITY SET}
\end{center}
\vskip0,5cm
\begin{center}
Ekaterina Vl. Bulinskaya\footnote{ \emph{Email address:} {\tt
bulinskaya@yandex.ru}}$^,$\footnote{The work is partially supported
by Dmitry Zimin Foundation ``Dynasty'' and RFBR grant 14-01-00318.}
\vskip0,2cm \emph{Lomonosov Moscow State University}
\end{center}
\vskip1cm

\begin{abstract}

We study a catalytic branching process (CBP) with any finite set of
catalysts. This model describes a system of particles where the
movement is governed by a Markov chain with arbitrary finite or
countable state space and the branching may only occur at the points
of catalysis. The results obtained generalize and strengthen those
known in cases of CBP with a single catalyst and of branching random
walk on ${\bf Z}^d$, $d\in{\bf N}$, with a finite number of sources
of particles generation. We propose to classify CBP with $N$
catalysts as supercritical, critical or subcritical according to the
value of the Perron root of a specified $N\times N$ matrix. Such
classification agrees with the moment analysis performed here for
local and total particles numbers. By introducing the criticality
set $C$ we also consider the influence of catalysts parameters on
the process behavior. The proof is based on construction of
auxiliary multi-type Bellman-Harris processes with the help of
hitting times under taboo and on application of multidimensional
renewal theorems.

\vskip0,5cm {\it Keywords and phrases}: catalytic branching process,
classification, hitting times under taboo, moment analysis,
multi-type Bellman-Harris process.

\vskip0,5cm 2010 {\it AMS classification}: 60J80, 60J27.

\end{abstract}

\section{Introduction}
We consider the model of \emph{catalytic branching process} (CBP)
with a finite number of catalysts. It describes a system of
particles moving in space and branching only in the presence of
catalysts. More exactly, let at the initial time $t=0$ there be a
single particle that moves on some finite or countable set $S$
according to a continuous-time Markov chain
$\eta=\{\eta(t),t\geq0\}$ generated by infinitesimal matrix
$Q=(q(x,y))_{x,y\in S}$. When this particle hits a finite set of
catalysts ${W=\{w_1,\ldots,w_N\}\subset S}$, say at the site $w_k$,
it spends there random time having the exponential distribution with
parameter $\beta_k>0$. Afterwards the particle either branches or
leaves the site $w_k$ with probabilities $\alpha_k$ and $1-\alpha_k$
($0\leq\alpha_k<1$), respectively. If the particle branches (at the
site $w_k$), it dies and just before the death produces a random
non-negative integer number $\xi_{k}$ of offsprings located at the
same site $w_k$. If the particle leaves $w_k$, it jumps to the site
$y\neq w_k$ with probability $-(1-\alpha_k)q(w_k,y)q(w_k,w_k)^{-1}$
and continues its movement governed by the Markov chain $\eta$. All
newly born particles are supposed to behave as independent copies of
their parent.

We assume that the Markov chain $\eta$ is irreducible and the matrix
$Q$ is conservative (i.e. $\sum\nolimits_{y\in S}{q(x,y)}=0$ where
$q(x,y)\geq0$ for $x\neq y$ and $q(x,x)\in(-\infty,0)$ for any $x\in
S$). Denote by $f_k(s):={\sf E}{s^{\xi_k}}$, $s\in[0,1]$, the
probability generating function of $\xi_k$, $k=1,\ldots,N$. We will
employ the standard assumption of existence of a finite derivative
$f_k'(1)$, that is the finiteness of ${\bf E}{\xi_k}$, for any
$k=1,\ldots,N$.

Study of CBP with a single catalyst was initiated even in the XX
century (see, e.g., \cite{ABY_98}). In this regard we also mention a
recent paper \cite{DR_13} where the main tool for the moment
analysis of the process was the spine technique, i.e. ``many-to-few
lemma'', and renewal theory. Note that the generalization to an
arbitrary finite set of catalysts is not straightforward, since
there exists a "competition" between catalysts. An important special
case of several catalysts (where $S={\bf Z}^{d}$, $d\in{\bf N}$, and
the Markov chain $\eta$ is symmetric space-homogeneous random walk
with finite variance of jump sizes) was examined in \cite{Y_MPS_13}.
There the sufficient conditions for exponential growth of the
particles numbers were obtained by analyzing spectral properties of
evolution operators. We propose another way permitting to establish
more precise results.

For branching processes the natural interesting problem is the
analysis of asymptotic behavior (as $t\to\infty$) of the local and
total size of population at time $t$ (for various branching
processes without catalysts see, e.g., \cite{Sew_74}). Let $\mu(t)$
stand for the total number of particles existing in CBP at time
$t\geq0$. In a similar way we define local numbers $\mu(t;y)$ as
quantities of particles located at separate points $y\in S$ at time
$t$. In this paper our aim is three-fold. Firstly, we introduce a
classification of CBP (with $N$ catalysts) treating it as {\it
supercritical, critical or subcritical} whenever one has, for the
Perron root $\rho$ of a certain $N\times N$ matrix, $\rho>1$,
$\rho=1$ or $\rho<1$, respectively. Moreover, the criticality set
$C$ revealing the influence of catalysts strength is introduced and
characterized. Secondly, we implement the moment analysis of the
local and total particles numbers to justify the naturalness of the
proposed classification (indeed, one will see that the asymptotic
behavior of the moments of any order is determined essentially by
the introduced class of CBP). Thirdly, we consider some particular
cases of CBP and thereby show that our study not only generalizes
the results in previous works but even refines them.

Our approach consists in involving hitting times under taboo (see,
e.g., \cite{B_SAM_12}, \cite{B_SPL_13} and \cite{Chung_60}, Ch.2,
Sec.11) and construction of auxiliary Bellman-Harris branching
processes with at most $N(N+1)+1$ types of particles. This approach
is inspired by \cite{VTY} where the branching random walk on ${\bf
Z}$ with a single catalyst was investigated by means of introducing
hitting times (without taboo) and a due two-type Bellman-Harris
process. Note also that in \cite{Carmona_Hu_2013}, for the study of
a discrete-time branching random walk on ${\bf Z}$ with multiple
catalysts, the authors considered an embedded multi-type
Galton-Watson branching process, resulting in "forgetting/erasing
the time spent between catalysts". The latter approach is fruitful
for classification of branching random walk since multi-type
Galton-Watson and multi-type Bellman-Harris processes have the same
supercritical, critical or subcritical regimes. However, for
subsequent study of branching random walks (moment analysis, limit
theorems etc.) multi-type Galton-Watson process is not sufficient
and thus Ph.Carmona and Y.Hu used in \cite{Carmona_Hu_2013} another
technique such as "many-to-few lemmas". Our first auxiliary process
is constructed in such a way that the study of CBP can be mainly
reduced to analysis of the Bellman-Harris process. The number of
particles in this process is chosen to guarantee its
indecomposability and cannot be less than $N$ and need not be
greater than $N(N+1)$. It is well-known (see, e.g., \cite{Sew_74},
Ch.4, Sec.5, 6 and 7) that an indecomposable multi-type
Bellman-Harris process is classified as supercritical, critical or
subcritical according to the value of the Perron root of the mean
matrix. This is the foundation for classification of CBP.
Furthermore, Lemma \ref{L:Perron_root} below gives even more
convenient classification with the help of the Perron root of a
specified irreducible $N\times N$ matrix. One more auxiliary
Bellman-Harris process is indispensable for study of total particles
numbers in CBP in case of transient Markov chain $\eta$. This
process is taken to be decomposable with a final type of particles.
The treatment of the auxiliary Bellman-Harris processes allows us
not only to classify CBP but also to derive a system of renewal
equations for the mean local and total particles numbers in CBP.
Afterwards, to implement the moment analysis of the local and total
particles numbers we use multidimensional renewal theorems
established in \cite{Crump_70} and \cite{Mode_68_2}.

This approach has advantages since there is elaborated theory of
multi-type Bellman-Harris processes and a vast majority of its
results can be applied to the auxiliary Bellman-Harris processes
leading to the new results for CBP. Mention in passing some recent
works on multi-type Bellman-Harris processes, see, e.g.,
\cite{Jones_11}, \cite{VT_TMIAN_13} and \cite{YY_09}.

Observe also that CBP can be considered as a Markov branching
process with at most countably many types of particles, since the
location of a particle can be associated with its type. Theory of
branching processes with countably many types of particles, despite
of its long history (see, e.g., \cite{Moy_67}), has not been
systematized until now in view of its complexity. Some new papers
such as \cite{BL_12}, \cite{BZ_09}, \cite{HLN_13} and
\cite{Sagitov_13} make an important contribution to this research
direction. However, the study there does not cover the results
presented in our paper. Among other investigations of models
describing particles movement and breeding we refer to recent papers
\cite{Biggins_12}, \cite{HuShi_09} and \cite{Lifshits_13}. Branching
random walks analyzed there are homogeneous in space whereas the
main feature of CBP is spatial non-homogeneity. Some differences in
behavior of a homogeneous branching random walk on ${\bf Z}$ and its
catalytic counterpart are discussed in \cite{Carmona_Hu_2013}.
Concluding the introduction we mention the close relation between a
catalytic branching random walk on ${\bf Z}$ and a super-Brownian
motion with a single point catalyst (see, e.g., \cite{DF_94} and
\cite{FLG_95}).

Now we describe the structure of the paper. In section \ref{s:ABHP}
we specify the first auxiliary Bellman-Harris process and propose a
classification for CBP. Section~\ref{s:SCS} is devoted to solution
of the following problem. Assume that we fix the Markov chain
generator $Q$ in CBP and vary ``intensities'' of catalysts, e.g.,
${m_k:={\sf E}{\xi_k}}$, $k=1,\ldots,N$. What is the set
$C\subset{\bf R}^N_+$ such that CBP is critical iff
${m=(m_1,\ldots,m_N)\in C}$? In particular, what is the proportion
of ``weak'' and ``powerful'' catalysts (possessing small or large
values of $m_i$, respectively) in critical regime? In section
\ref{s:SCS} we obtain a complete description of $C$ by means of
equation involving determinants of some matrices, indicate the
smallest parallelepiped $[0,M_1]\times \ldots \times [0,M_N]$
containing $C$ and illustrate our approach by two pictures of the
set $C$ when $N=2$ and $N=3$. Observe that these plots are not a
product of simulation and display the result of direct computation.
Section \ref{s:MACBP} contains Theorem \ref{T:MACBP} and its proof
representing the moment analysis of the local and total particles
numbers in CBP. There we construct the second auxiliary
Bellman-Harris process required for study of the total size of
population whenever Markov chain $\eta$ is transient.
Section~\ref{s:Applications} demonstrates applications of our
results to the catalytic branching random walk on ${\bf Z}^d$ and to
the branching process with a single catalyst. A detailed comparison
of our results with those known before is given as well.

\section{Auxiliary Bellman-Harris process}\label{s:ABHP}
Let us briefly describe a Bell\-man-Harris branching process with
particles of $n$ types, $n\in{\bf N}$. It is initiated by a single
particle of type $i = 1,\ldots,n$. The parent particle has a random
life-length with a cumulative distribution function (c.d.f.)
$G_i(t)$, $t\geq0$. When dying the particle produces offsprings
according to a probability generating function $g_i({\bf s})$, ${\bf
s}=(s_1,\ldots,s_n)\in[0,1]^n$. The new particles of type
${j=1,\ldots,n}$ evolve independently with the life-length
distribution $G_j(\cdot)$ and an offspring generating function
$g_j(\cdot)$. Let
$$M:=\left(m_{i,j}\right)_{i,j=1}^n\;\; \mbox{with}\;\;
m_{i,j}=\left.\partial_{s_{j}}g_{i}\right|_{{\bf s}=(1,\ldots,1)}
$$
be the {\it mean matrix} of the process. The Bellman-Harris
branching process is called {\it indecomposable} if the non-negative
matrix $M$ is {\it irreducible} (for the latter notion, see, e.g.,
\cite{Seneta_06}, Ch.1, Sec.3). Moreover, if its Perron root $\rho$
(i.e. eigenvalue having the maximal modulus) is such that  $\rho>1$,
$\rho=1$ or $\rho<1$, then the Bellman-Harris process is called
supercritical, critical or subcritical, respectively (see, e.g.,
\cite{Sew_74}, Ch.4, Sec.5, 6 and 7). Denote the number of particles
of type $j$ existing at time $t$ by $Z_j(t)$, $t\geq0$,
$j=1,\ldots,n$.

Before demonstrating how an auxiliary Bellman-Harris process can be
constructed in the framework of CBP we have to introduce some
notation. Consider a particle moving on the set $S$ in accordance
with the Markov chain generated by $Q$ and starting at state $x$.
Let $_H\overline{\tau}_{x,y}$, $x,y\in S$, $H\subset S$, be the time
spent by the particle after leaving the starting point $x$ until the
first hitting $y$ if the particle's trajectory does not pass $H$.
Otherwise (if the particle's trajectory passes $H$ before the first
hitting $y$), ${_H\overline{\tau}_{x,y}=\infty}$. The (extended)
random variable $_H\overline{\tau}_{x,y}$ is called a {\it hitting
time} of state $y$ {\it under taboo} $H$ after the first exit out of
the starting state $x$. Denote by $_H\overline{F}_{x,y}(t)$,
$t\geq0$, the improper c.d.f. of $_H\overline{\tau}_{x,y}$.
Evidently,
${_H\overline{F}_{x,y}(0)=(\delta_{x,y}-1)q(x,y)q(x,x)^{-1}}$ where
$\delta_{x,y}$ is the Kronecker delta. Explicit formulae for the
probability of finiteness of $_H\overline{\tau}_{x,y}$, i.e. for
$_H\overline{F}_{x,y}(\infty)=\lim\nolimits_{t\to\infty}{_H\overline{F}_{x,y}(t)}$,
via taboo probabilities and Green's function were derived in
\cite{B_SPL_13}. Whenever the taboo set $H$ is empty we write
$\overline{\tau}_{x,y}$ and $\overline{F}_{x,y}(\cdot)$ instead of
${_\emptyset\overline{\tau}_{x,y}}$ and
${_\emptyset\overline{F}_{x,y}(\cdot)}$, respectively.

Return to CBP. We tentatively assume that CBP starts at $w_{p}$ for
some ${p=1,\ldots,N}$. Set
$K_j:=\{k=1,\ldots,N:{_{W_k}\overline{F}}_{w_j,w_k}(\infty)-{_{W_k}\overline{F}_{w_j,w_k}}(0)>0\}$
for $W_k:=W\setminus\{w_k\}$ and $j=1,\ldots,N$. Let
$K_j=\{k(1,j),\ldots,k(|K_j|,j)\}$ where $1\leq
k(1,j)<k(2,j)<\ldots<k(|K_j|,j)\leq N$ and $|\cdot|$ is the
cardinality of a finite set. We divide the particles population
existing at time $t\geq0$ into $L+1$ groups with
$L:=N+\sum\nolimits_{j=1}^{N}{|K_{j}|}$. The particles located at
time $t$ at $w_{j}$ form the $j$-th group having cardinality
$\mu(t;w_{j})$, $j=1,\ldots,N$. Consider a family consisting of
particles which have left $w_{j}$ at least once within time interval
$[0,t]$, upon the last leaving $w_{j}$ have not yet reached $W$ by
time $t$ but eventually will hit $w_{k(i,j)}$ before possible
hitting $W_{k(i,j)}$, $i=1,\ldots,|K_j|$, $j=1,\ldots,N$. This
family has cardinality denoted by $\mu_{j,i}(t)$ and corresponds to
the group number $(L(j)+i)$ where
$L(j):=N+\sum\nolimits_{l=1}^{j-1}{|K_{l}|}$. The group number $L+1$
comprises the rest of particles not included into the above $L$
groups. Note that the last group consists of the particles having
infinite life-lengths since after time $t$ they will not hit the set
of catalysts any more. So, after time $t$ these particles will not
produce any offsprings and have no influence on the numbers of
particles in other $L$ groups.

Now we can introduce an auxiliary Bellman-Harris process to employ
it for the study of CBP. Consider an $L$-dimensional Bellman-Harris
process starting with a single particle of type $p$ and having the
following c.d.f. and offspring generating functions
\begin{equation}\label{Gg_definition_1}
G_j(t)=1-e^{-\beta_jt},\quad g_j({\bf
s})=\alpha_j
f_j(s_j)+(1-\alpha_j)\sum\limits_{k=1}^{N}{_{W_k}\overline{F}_{w_j,w_k}(0)s_k}
\end{equation}
$$
+(1-\alpha_j)\!\sum\limits_{i=1}^{|K_j|}T_{i,j}(\infty)s_{L(j)+i}
+(1-\alpha_j)\left(1
-\sum\limits_{k=1}^{N}{_{W_k}\overline{F}_{w_j,w_k}(\infty)}\right),
$$
\begin{equation}\label{Gg_definition_2}
G_{L(j)+i}(t)=\frac{T_{i,j}(t)}{T_{i,j}(\infty)},\quad
g_{L(j)+i}({\bf s})=s_{k(i,j)}
\end{equation}
where $T_{i,j}(t):={_{W_{k(i,j)}}
\overline{F}_{w_j,w_{k(i,j)}}(t)}-{_{W_{k(i,j)}}\overline{F}_{w_j,w_{k(i,j)}}(0)}$,
$0\leq t\leq \infty$ and $i=1,\ldots,|K_j|$, $j=1,\ldots,N$. It is
easy to see that for the process constructed in this way the vector
$(Z_1(t),\ldots,Z_L(t))$, at each $t\geq0$, has the same
distribution as the vector whose $j$-th component is $\mu(t;w_j)$ and
$(L(j)+i)$-th component is $\mu_{j,i}(t)$, $i=1,\ldots,|K_j|$,
$j=1,\ldots,N$.

Mean matrix $M=(m_{k,l})_{k,l=1}^L$ of the introduced Bellman-Harris
process is of block form, i.e.
\begin{equation}\label{matrix_M=}
M=\left(
      \begin{array}{cccc}
        M_{1,1} & M_{1,2} & \ldots & M_{1,N+1} \\
        M_{2,1} & M_{2,2} & \ldots & M_{2,N+1} \\
        \ldots & \ldots & \ldots & \ldots \\
        M_{N+1,1} & M_{N+1,2} & \ldots & M_{N+1,N+1} \\
      \end{array}
    \right).
\end{equation}
Here the $N\times N$ matrix $M_{1,1}=(m_{k,l})_{k,l=1}^N$ possesses
the following entries $m_{k,l}=\delta_{k,l}\,\alpha_k
f'_k(1)+(1-\alpha_k){_{W_l}\overline{F}_{w_k,w_l}}(0)$. For
$j=1,\ldots,N$, the elements of $N\times|K_j|$ matrix $M_{1,j+1}$
vanish everywhere except for the $j$-th row and for
${i=1,\ldots,|K_j|}$ one has
$m_{j,L(j)+i}=(1-\alpha_j)T_{i,j}(\infty)$. When $i=1,\ldots,N$ the
$|K_i|\times N$ matrix $M_{i+1,1}$ is such that omitting its $j$-th
null columns for all $j\notin K_i$ one gets the $|K_i|\times|K_i|$
identity matrix. For $i,j=1,\ldots,N$, the $|K_i|\times|K_j|$ matrices
$M_{i+1,j+1}$ have zero entries.

Let us verify that the specified Bellman-Harris process is
indecomposable. In view of irreducibility of the Markov chain $\eta$
there exists a finite path from $w_k$ to $w_l$ with a positive
probability for each $k,l=1,\ldots,N$. In the framework of CBP such
path has also positive probability, since being at catalyst site
$w_i$ a particle can leave it without branching with positive
probability $1-\alpha_i$ for each $i=1,\ldots,N$. Among the sites
visited successively by that path of $\eta$ choose those from $W$,
say, $w_{p(0)},w_{p(1)},\ldots,w_{p(J)},w_{p(J+1)}$ with $p(0)\!:=k$
and ${p(J+1)\!:=l}$. By the construction of auxiliary Bellman-Harris
process this path corresponds to transformations of a particle type
from $k$ to $l$. Namely, if the path of $\eta$ hits $w_{p(j+1)}$
immediately after leaving $w_{p(j)}$ then the particle type in
Bellman-Harris process changes from $p(j)$ to $p(j+1)$. Otherwise,
$p(j+1)\in K_{p(j)}$, i.e. $p(j+1)=k(i,p(j))$ for some
$i=1,\ldots,|K_{p(j)}|$, and the particle type change from $p(j)$ to
$p(j+1)$ involves the intermediate type $L(p(j))+i$. Hence, we
conclude that, for each $k$ and $l$ from $\{1,\ldots,N\}$, there
exists $n=n(k,l)\in{\bf N}$ and a collection
$\{r(1),\ldots,r(n)\}\subset\{1,\ldots,L\}$ such that
$m_{k,r(1)}m_{r(1),r(2)}\ldots m_{r(n),l}>0$. Since by matrix $M$
definition in (\ref{matrix_M=}) one has $m_{L(j)+i,k(i,j)}=1>0$ and
$m_{j,L(j)+i}>0$ for each $i=1,\ldots,|K_j|$ and $j=1,\ldots,N$, the
previous statement holds true for $k,l=1,\ldots,L$ as well.
Consequently, we have checked that for each $k$ and $l$ from
$\{1,\ldots,L\}$, there exists $n=n(k,l)\in{\bf N}$ such that
$m^{(n)}_{k,l}>0$ where $m^{(n)}_{k,l}$ is the $(k,l)$-th element of
$M^n$. So, we have verified that the matrix $M$ is irreducible.
Moreover, if we consider the nontrivial case $f_i'(1)>0$ for some
$i=1,\ldots,N$, then $m_{i,i}>0$ and the matrix $M$ is
\emph{acyclic} (or \emph{aperiodic}) (see, e.g., \cite{Seneta_06},
Ch.1, Sec.2). Thereby, in this case we have shown that $M$ is a
\emph{primitive} matrix (see, e.g., \cite{Seneta_06}, Ch.1, Sec.3).

On account of the Perron-Frobenius theorem for irreducible matrices
(see, e.g., \cite{Seneta_06}, Ch.1, Sec.4) the matrix $M$ has a
positive real eigenvalue $\rho(M)$ of maximal modulus which is
called the \emph{Perron root} of $M$. Using the classification of
the constructed Bellman-Harris process we call CBP
\emph{supercritical}, \emph{critical} or \emph{subcritical} if
$\rho(M)>1$, $\rho(M)=1$ or $\rho(M)<1$, respectively, where $M$ is
specified by (\ref{matrix_M=}). Since $M$ is independent of the CBP
starting point, the same is true for the classification of CBP. So,
from here on we omit our assumption that CBP starts at some $w_p\in
W$.

Let $G^\ast_i(\lambda):=\int\nolimits_{0-}^{\infty}{e^{-\lambda
t}\,dG_i(t)}$, $\lambda\geq0$, be the Laplace-Stieltjes transform of
c.d.f. $G_i$ and
$H(\lambda):=\left(G^\ast_i(\lambda)m_{i,j}\right)_{i,j=1}^{L}$.
Note that $H(0)=M$. Matrix $H(\lambda)$ is irreducible in view of
irreducibility of $M$. Put also
$D(\lambda)=(d_{i,j}(\lambda))_{i,j=1}^N$, $\lambda\geq0$, with
$d_{i,j}(\lambda):=\delta_{i,j}\,\alpha_i
f_i'(1)G^\ast_i(\lambda)+(1-\alpha_i)G^\ast_i(\lambda){_{W_j}\overline{F}^\ast_{w_i,w_j}(\lambda)}$.
It is easily seen that the matrix $D=(d_{i,j})_{i,j=1}^N$ defined as
$D:=D(0)$ has entries
\begin{equation}\label{entries_D}
{d_{i,j}=\delta_{i,j}\,\alpha_i
f_i'(1)+(1-\alpha_i){_{W_j}\overline{F}_{w_i,w_j}(\infty)}}.
\end{equation}
For an irreducible matrix $A$, let $\rho(A)$ stand for the Perron
root of $A$. The following statement gives a convenient criticality
condition for CBP.
\begin{Lm}\label{L:Perron_root}
For each $\lambda\geq0$, the matrix $D(\lambda)$ is irreducible and
the values $\rho(D(\lambda))$ and $\rho(H(\lambda))$ are
simultaneously greater than $1$, equal to $1$ or less than $1$. In
particular, the Perron roots $\rho(D)$ and $\rho(M)$ are either both
greater than $1$, or equal to $1$, or are less than $1$.
\end{Lm}
{\sc Proof.} The irreducibility of $D$ (and, therefore, of
$D(\lambda)$) is established in the same manner as for matrix $M$.
In accordance with the Perron-Frobenius theorem matrix $H(\lambda)$
has a strictly positive left eigenvector ${\bf v}$ corresponding to
the eigenvalue $\rho(H(\lambda))$. Consequently,
$H(\lambda)^\top{\bf v}=\rho(H(\lambda)){\bf v}$ where $\top$ means
transposition and ${\bf v}$ is considered as a column vector. In
other words, $\left(H(\lambda)^\top-\rho(H(\lambda))I\right){\bf
v}={\bf 0}$ for $I$ being the identity matrix and ${\bf 0}$ standing
for the zero of ${\bf R}^L$. Apply some equivalent transformations
to the obtained system of equations, namely, for each
${j,k=1,\ldots,N}$ add the $(L(j)+i)\mbox{-th}$ row multiplied by
$G^\ast_{L(j)+i}(\lambda)\rho(H(\lambda))^{-1}$ to the $k$-th row if
$k\in K_j$, i.e. ${k=k(i,j)}$ for some $i=1,\ldots,|K_j|$. Focusing
on the first $N$ equations only, we deduce that
$${\bf v}_0^\top\left(H_{1,1}(\lambda)+\rho(H(\lambda))^{-1}R(\lambda)-\rho(H(\lambda))I\right)={\bf
0}^\top.$$ Here we denote
$H_{1,1}(\lambda):=(G^\ast_i(\lambda)m_{i,j})_{i,j=1}^N$ and matrix
$R(\lambda):=(r_{i,j}(\lambda))_{i,j=1}^N$ has the following entries
${r_{i,j}(\lambda):=(1-\alpha_i)G^\ast_i(\lambda)\left({_{W_j}}\overline{F}^\ast_{w_i,w_j}(\lambda)-
{_{W_j}}\overline{F}_{w_i,w_j}(0)\right)}$ whereas ${\bf v}_0$ is a
strict\-ly positive vector consisting of the first $N$ coordinates
of ${\bf v}$ and ${\bf 0}$ means the zero of ${\bf R}^N$. Note that
$D(\lambda)=H_{1,1}(\lambda)+R(\lambda)$. Evidently, if
$\rho=\rho(H(\lambda))\geq1$ then
$$\rho^{-1}{\bf v}_0^\top\left(D(\lambda)-\rho^2
I\right)={\bf
v}_0^\top\left(\rho^{-1}H_{1,1}(\lambda)+\rho^{-1}R(\lambda) -\rho
I\right)$$
$$\leq {\bf v}_0^\top\left(H_{1,1}(\lambda)\!+\!\rho^{-1}R(\lambda)\!-\!\rho I\right)\leq
{\bf v}_0^\top\left(H_{1,1}(\lambda)\!+\!R(\lambda)\!-\!\rho
I\right)={\bf v}_0^\top\left(D(\lambda)\!-\!\rho I\right).$$ In a
similar way, if $\rho(H(\lambda))\leq1$ then the previous chain of
inequalities holds true when each symbol $\leq$ is replaced by
$\geq$. Finally, we get
\begin{equation}\label{>=}
D(\lambda)^\top {\bf v}_0\geq\rho(H(\lambda)){\bf
v}_0\;\:\mbox{and}\;\: D(\lambda)^\top {\bf
v}_0\leq\rho(H(\lambda))^2{\bf v}_0\;\:
\mbox{for}\;\:\rho(H(\lambda))\geq1,
\end{equation}
\begin{equation}\label{<=}
D(\lambda)^\top {\bf v}_0\leq\rho(H(\lambda)){\bf
v}_0\;\:\mbox{and}\;\: D(\lambda)^\top {\bf
v}_0\geq\rho(H(\lambda))^2{\bf
v}_0\;\:\mbox{for}\;\:\rho(H(\lambda))\leq1.
\end{equation}
On account of Theorem 1.6 in \cite{Seneta_06}, Ch.1, Sec.4, the
second relation in (\ref{>=}) and the first one in (\ref{<=}) imply
that ${\rho(D(\lambda))\leq\rho(H(\lambda))^2}$ for
${\rho(H(\lambda))\geq1}$ and
${\rho(D(\lambda))\leq\rho(H(\lambda))}$ for
$\rho(H(\lambda))\leq1$. Examining the proof of the Perron-Frobenius
theorem we conclude that due to the first inequality in (\ref{>=})
and the second  in (\ref{<=}) one has
$\rho(D(\lambda))\geq\rho(H(\lambda))$ for $\rho(H(\lambda))\geq1$
and correspondingly ${\rho(D(\lambda))\geq\rho(H(\lambda))^2}$ for
$\rho(H(\lambda))\leq1$. Consequently, if ${\rho(H(\lambda))\geq1}$
then $\rho(H(\lambda))\leq\rho(D(\lambda))\leq\rho(H(\lambda))^2$
whereas if $\rho(H(\lambda))\leq1$ we obtain
${\rho(H(\lambda))^2\leq\rho(D(\lambda))\leq\rho(H(\lambda))}$.
These estimates entail the assertion of Lemma~\ref{L:Perron_root}.
$\square$

\begin{Def}
CBP is called supercritical, critical or subcritical whenever
$\rho(D)>1$, $\rho(D)=1$ or $\rho(D)<1$, respectively, where $D$ is
the matrix  specified by way of $(\ref{entries_D})$.
\end{Def}

\section{Structure of the criticality set}\label{s:SCS}
In this section we focus on the nonsingular case $\alpha_i\in(0,1)$
for each $i=1,\ldots,N$. One can see that, given a Markov chain
generator $Q$ and a collection $\alpha_1,\ldots,\alpha_N$, the
matrix $D$ depends on $m_{i}={\bf E}{\xi_i}=f_i'(1)$ only and not on
the explicit form of $f_i(\cdot)$, $i=1,\ldots,N$, that is
$D=D(m_1,\ldots,m_N)$. Define \emph{the criticality set} as
$${C:=\{(m_1,\ldots,m_N)\in{\bf
R}_+^N:\rho(D)=1\}}.$$ Study of the criticality set structure allows
us to understand the relationship between ``weak'' and ``powerful''
catalysts (possessing small or large values of $m_i$) leading to the
CBP  supercriticality, criticality or subcriticality.

Clearly, $C\subset\{(m_1,\ldots,m_N)\in{\bf
R}_+^N:\det\left(D-I\right)=0\}$ though the latter set also covers
the case when $N>1$, $\rho(D)>1$ and $D$ has an eigenvalue $1$.
Using the Laplace expansion of determinants we derive that for $N>1$
one has $\det\left(D-I\right)=0$ iff
$$\left(\alpha_i
m_i-1\right)\det(D-I)_{i,i}=
(1-\alpha_i)\sum\limits_{j=1}^{N}{(-1)^{i+j-1}{_{W_j}\overline{F}_{w_i,w_j}}(\infty)
\det\left(D-I\right)_{i,j}}$$ where $\left(D-I\right)_{i,j}$ is the
$(N-1)\times(N-1)$ matrix that results from deleting the $i$-th row
and the $j$-th column of $D-I$, $i,j=1,\ldots,N$. Since matrix $D$ is
irreducible by Lemma~\ref{L:Perron_root}, the maximal eigenvalue of
any principal sub-matrix of $D$ is strictly less than $\rho(D)$
(see, e.g., \cite{Gantmacher_00}, vol.2, Ch.13, Sec.3). Thus, if
$\rho(D)=1$ then $\det(D-I)_{i,i}\neq0$ for any $i=1,\ldots,N$ and
\begin{eqnarray}
m_i&=&\frac{1-(1-\alpha_i){_{W_i}\overline{F}_{w_i,w_i}}(\infty)}{\alpha_i}\label{m_i=det}\\
&+&\frac{(1-\alpha_i)\sum\nolimits_{j=1,\,j\neq
i}^{N}{(-1)^{i+j-1}{_{W_j}\overline{F}_{w_i,w_j}}(\infty)
\det\left(D-I\right)_{i,j}}}{\alpha_i\det(D-I)_{i,i}}.\nonumber
\end{eqnarray}
Obviously, the matrix $\left(D(m_1,\ldots,m_N)-I\right)_{i,j}$ does
not depend on $m_i$ and $m_j$ for any $i,j=1,\ldots,N$. Therefore,
if $(m_1,\ldots,m_n)\in C$ and all $m_j$, for ${j=1,\ldots,N}$,
$j\neq i$, are fixed then $m_i$ is uniquely determined by
(\ref{m_i=det}).

Let us find the bounds for possible values of $m_i$, $i=1,\ldots,N$,
such that $(m_1,\ldots,m_N)\in C$. At first recall that according to
the Perron-Frobenius theorem, for irreducible matrices $A$ and $B$,
the Perron root $\rho(A)$ does not exceed the Perron root $\rho(B)$
whenever all the elements of matrix $B-A$ are nonnegative.
Furthermore, then $A=B$ if ${\rho(A)=\rho(B)}$. Assume that
$m_{i}=0$ for all $i=1,\ldots,N$. In this case the maximal row sum
of $D$ is less than 1 and, consequently, $\rho(D)<1$ in view of
Theorem 1.5 and Corollary 1 in \cite{Seneta_06}, Ch.1, Sec.1 and~4.
If now we let $m_1\to\infty$ then the Perron root
$\rho(D(m_1,0,\ldots,0))$ strictly grows to infinity by virtue of
the Perron-Frobenius theorem proof. Since $\rho(D)$ is continuously
dependent on any element of matrix $D$, there exists $M_1>0$ such
that ${\rho\left(D(M_1,0,\ldots,0)\right)=1}$. Moreover, if
$m_1<M_1$ then $\rho\left(D(m_1,0,\ldots,0)\right)<1$ whereas for
$m_1>M_1$ and any nonnegative $m_2,\ldots,m_N$ one can verify that
$\rho\left(D(m_1,m_2,\ldots,m_N)\right)>1$. The exact value of $M_1$
can be found  by setting $m_2=m_3=\ldots=m_N=0$ in (\ref{m_i=det}).
In a similar way, we may indicate positive $M_2,\ldots,M_N$ such
that $(m_1,\ldots,m_N)\in C$ with $m_i=M_i$ and $m_j=0$ for each
$i=2,\ldots,N$ and $j\neq i$. Thus, we have found the smallest
parallelepiped $[0,M_1]\times\ldots\times[0,M_N]$ containing $C$,
i.e. $C\subset[0,M_1]\times\ldots\times[0,M_N]$.

Assume now that $m_1=m_1^0$ where $m_1^0\in[0,M_1)$. Reasoning as
above we conclude that there exists $M_2(m_1^0)\in(0,M_2]$ such that
the equality ${\rho\left(D(m_1^0,M_2(m_1^0),0,\ldots,0)\right)=1}$
holds true. The value $M_2(m_1^0)$ is given by formula
(\ref{m_i=det}) when $m_1=m_1^0$ and $m_3=\ldots=m_N=0$. Let
${m_2=m_2^0}$ with $m_2^0\in[0,M_2(m_1^0))$. Next we can take
$M_3(m_1^0,m_2^0)\in(0,M_3]$ so that
$\rho\left(D(m_1^0,m_2^0,M_3(m_1^0,m_2^0),0,\ldots,0))\right)=1$.
Using identity (\ref{m_i=det}) and by way of the above reasoning we
choose the array $m_i^0\in[0,M_i(m_1^0,\ldots,m_{i-1}^0))$ and thus
reveal a collection of numbers
$M_i(m_1^0,\ldots,m_{i-1}^0)\in(0,M_i]$ such that
${\rho\left(D(m_1^0,\ldots,m_{i-1}^0,M_i(m_1^0,\ldots,m_{i-1}^0),0,\ldots,0)\right)=1}$
for $i=1,\ldots,N-1$ (naturally, we set $M_1(\emptyset)=M_1$).
Finally, to pick $m_N^0$ for ${(m_1^0,\ldots,m_N^0)\in C}$ we employ
(\ref{m_i=det}) with $m_1=m_1^0$, $\ldots$, $m_{N-1}=m_{N-1}^0$.
Thereby we demonstrate how one can vary a collection
$(m_1^0,\ldots,m_N^0)$ providing the criticality of the
corresponding CBP. Note that if for some step ${i=1,\ldots,N-1}$ we
choose $m_i^0= M_i(m_1^0,\ldots,m_{i-1}^0)$ then $m_k^0=0$ for all
$k=i+1,\ldots,N$. Moreover, if for some step $i=1,\ldots,N-1$ we
take ${m_i>M_i(m_1^0,\ldots,m_{i-1}^0)}$ then CBP is supercritical
for any nonnegative values $m_{i+1},\ldots,m_N$. At last, if
$m_N^0>0$ then the choice $(m_1^0,\ldots,m_{N-1}^0,m_N)$ with
$m_N>m_N^0$ or $m_N<m_N^0$ leads to supercritical or subcritical
CBP, respectively.

\begin{figure}
\includegraphics[width=12.5cm]{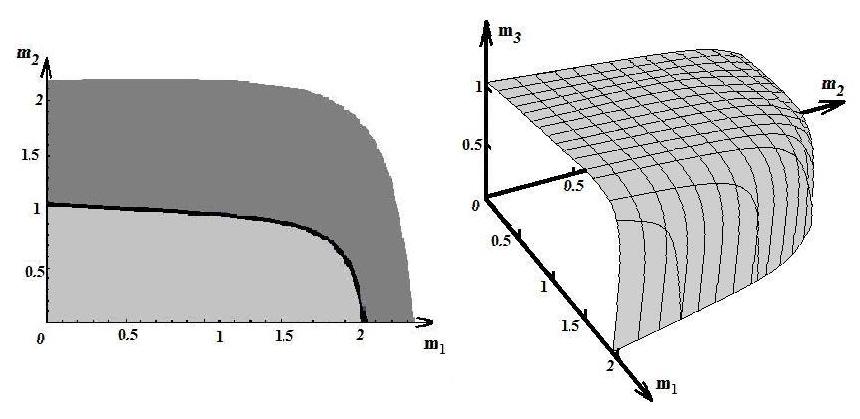}
\caption[]{Examples of the criticality set $C$ for $N=2$ and
$N=3$.}\label{criticality_set}
\end{figure}

To illustrate this discussion we provide below Figure
\ref{criticality_set} depicting the set $C$ for some particular
cases. On the first picture there is a black curve representing the
set $C$ for a simple random walk on ${\bf Z}$ with two catalysts
located at neighboring points. In this case, for instance, by
Theorem~2 in \cite{B_SAM_12}, we deduce that
${_{W_j}\overline{F}_{w_i,w_j}}(\infty)=1/2$ for each $i,j=1,2$. We
also set $\alpha_1=0.3$ and $\alpha_2=0.8$. The domain marked by
light grey color under the curve is related to subcritical CBP
whereas the domain marked by dark grey color over the curve
corresponds to supercritical CBP. On the second picture there is a
surface depicting the set $C$ for a simple random walk on ${\bf Z}$
with three catalysts located at three subsequent points. With the
help of Theorem~2 in \cite{B_SPL_13} we derive that
${_{W_j}\overline{F}_{w_i,w_j}}(\infty)$ equals either $1/2$, given
that $|w_i-w_j|\leq1$, or $0$, provided that $|w_i-w_j|=2$,
$i,j=1,2,3$, except for the case
${_{W_2}\overline{F}_{w_2,w_2}}(\infty)=0$. We also assume that
$\alpha_1=0.3$, $\alpha_2=0.6$ and $\alpha_3=0.8$. The domain under
the surface corresponds to subcritical CBP while that over the
surface is related to supercritical CBP (the net on the surface
displays the level lines).

\section{Moment analysis of CBP}\label{s:MACBP} It was Ch.J.~Mode
who showed (see \cite{Mode_68_1}) that the mean particles numbers in
a supercritical indecomposable multi-type Bellman-Harris process
grow exponentially (as time tends to infinity) with a certain rate
$\nu$ called a Malthusian parameter. The parameter $\nu$ is positive
and determined as the unique solution to equation $\rho(H(\nu))=1$.
Lemma \ref{L:Perron_root} demonstrates that for our auxiliary
$L$-dimensional supercritical Bellman-Harris process this procedure
can be simplified and the Malthusian parameter can be found as the
unique solution to equation $\rho(D(\nu))=1$ (observe that square
matrix $H(\lambda)$ is of size $L\geq N$ whereas square matrix
$D(\lambda)$ is of size $N$). So, in the framework of CBP the
Malthusian parameter $\nu$ is defined as the unique (positive)
solution to equation $\rho(D(\nu))=1$ whenever $\rho(D)>1$. Given
$\rho(D)=\rho(D(0))>1$, such solution always exists, since according
to the Perron-Frobenius theorem the function $\rho(D(\lambda))$ in
variable $\lambda$ is strictly decreasing and tends to $0$ as
$\lambda\to\infty$.

Recall that, for an $n\times n$ matrix $R$, the matrix $R_{i,j}$ or
$(R)_{i,j}$ is the ${(n-1)\times(n-1)}$ matrix that results from
deleting the $i$-th row and the $j$-th column of $R$,
$i,j=1,\ldots,n$. To formulate our main result we need some new
notation. For supercritical CBP consider functions $a_n(x,y)$,
$x,y\!\in\! S$, $n\!\in\!\bf N$, which appear in description of the
asymptotic behavior, as $t\to\infty$, of the $n$-th factorial
moments ${m_n(t;x,y):={\sf
E}_x\mu(t;y)(\mu(t;y)-1)\ldots(\mu(t;y)-n+1)}$ of the local
particles numbers in CBP (index $x$ denotes the starting point of
CBP). The function $a_1(x,y)$ has the most concise form when $x,y\in
W$, namely,
\begin{equation}\label{a_1(w_i,w_j)_definition}
a_1(w_i,w_j)=\frac{\Delta_{j,i}(\nu)}{(\nu+\beta_j)
\Delta'(\nu)},\quad i,j=1,\ldots,N,
\end{equation}
where
$\Delta_{i,j}(\lambda)\!:=\!(-1)^{i+j}\det(I-D(\lambda))_{i,j}$ and
$\Delta'(\nu)\!:=\!\left.\frac{d}{d\lambda}\det(I-D(\lambda))\right|_{\lambda=\nu}$.
If $N=1$ then $\Delta_{1,1}(\lambda)$, $\lambda\geq0$, is supposed
to be identically $1$. For $x\in W$ and $y\in S\setminus W$ we set
\begin{equation}\label{a_1(w_i,y)_definition}
a_1(w_i,y)=\frac{\sum_{j=1}^N\Delta_{j,i}(\nu)(1-\alpha_j)G^\ast_j(\nu){_W\overline{F}^\ast_{w_j,y}(\nu)}}
{(\nu-q(y,y))\Delta'(\nu)(1-{_W F^\ast_{y,y}(\nu)})},\quad
i=1,\ldots,N.
\end{equation}
Here ${_H
F^\ast_{x',y'}(\lambda)}:=(-q(x',x'))(\lambda-q(x',x'))^{-1}{_H\overline{F}^\ast_{x',y'}(\lambda)}$,
$\lambda\geq0$, $x',y'\in S$, $H\subset S$, i.e. ${_H
F_{x',y'}(t)}$, $t\geq0$, is the c.d.f. of a hitting time of state
$y'$ under taboo $H$ given that the Markov chain $\eta$ starts at
state $x'$ (see, e.g., \cite{B_SPL_13}). Note that
$G^\ast_j(\lambda)=\beta_j(\lambda+\beta_j)^{-1}$ for each
$\lambda\geq0$ and $j=1,\ldots,N$. The values $a_n(x,y)$ for $n>1$,
$x\in W$ and $y\in S$ are evaluated with the help of the following iterative
scheme
\begin{equation}\label{a_n(w_i,y)_definition}
a_n(w_i,y)=\sum\limits_{k=1}^N\frac{\beta_k\Delta_{k,i}(n\nu)}{(n\nu+\beta_k)\Delta(n\nu)}
h_{n,k}(a_1(w_k,y),\ldots,a_{n-1}(w_k,y))
\end{equation}
where $i=1,\ldots,N$, $\Delta(\lambda):=\det(I-D(\lambda))$,
$\lambda\geq0$,
\begin{equation}\label{h_n,k_definition}
h_{n,k}(z_1,\ldots,z_{n-1}):=\alpha_k\sum\limits_{r=2}^{n}\frac{f^{(r)}_k(1)}{r!}
\sum\limits_{\substack{i_1,\ldots,i_r>0,\\
i_1+\ldots+i_r=n}}\frac{n!}{i_1!\ldots i_r!}z_{i_1}\ldots z_{i_r}
\end{equation}
and $z_1,\ldots,z_{n-1}\geq0$, $k=1,\ldots,N$. In this way the
values $a_n(x,y)$ for $n\geq1$, $x\in S\setminus W$ and $y\in S$ are
determined by formula
\begin{equation}\label{a_n(x,y)_definition}
a_n(x,y)=\sum\nolimits_{i=1}^N{_{W_i}F^\ast_{x,w_i}(n\nu)}a_n(w_i,y).
\end{equation}

Define functions $A_n(x)$, $x\in S$, $n\in{\bf N}$, arising in the
asymptotic behavior of the $n$-th factorial moments ${M_n(t;x):={\sf
E}_x\mu(t)(\mu(t)-1)\ldots(\mu(t)-n+1)}$ of the total particles
numbers in supercritical CBP. Namely, for $n=1$ and $x\in W$, put
\begin{equation}\label{A_1(w_i)_definition}
A_1(w_i)=\sum\limits_{j=1}^N\frac{\alpha_j\beta_j(f'_j(1)-1)\Delta_{j,i}(\nu)}
{\nu(\nu+\beta_j)\Delta'(\nu)},\quad i=1,\ldots,N.
\end{equation}
When $n>1$ and $x\in W$ the values $A_n(x)$ are evaluated according to
the following iterative scheme
\begin{equation}\label{A_n(w_i)_definition}
A_n(w_i)=\sum_{j=1}^N\frac{\beta_j\Delta_{j,i}(n\nu)}{(n\nu+\beta_j)\Delta(n\nu)}
h_{n,j}(A_1(w_j),\ldots,A_{n-1}(w_j))
\end{equation}
where $i=1,\ldots,N$. Bearing on the latter equalities the values
$A_n(x)$ for $n\geq1$ and $x\in S\setminus W$ are obtained  by
way of
\begin{equation}\label{A_n(x)_definition}
A_n(x)=\sum\nolimits_{i=1}^N{_{W_i}F^\ast_{x,w_i}(n\nu)}A_n(w_i).
\end{equation}

For critical CBP introduce functions $b_n(x,y)$ and $B_n(x)$,
$x,y\in S$, $n\in{\bf N}$, appearing in the asymptotic behavior of
$m_n(t;x,y)$ and $M_n(t;x)$, respectively. Namely, for
$i,j=1,\ldots,N$ and $y\in S\setminus W$, set
\begin{equation}\label{b_1(w_i,w_j)_definition}
b_1(w_i,w_j)=\frac{\Delta_{j,i}(0)}{\beta_j\Delta'(0)},\quad
b_1(w_i,y)=\frac{\sum_{j=1}^N\Delta_{j,i}(0)(1-\alpha_j){_WF_{w_j,y}(\infty)}}
{(-q(y,y))\Delta'(0)(1-{_W F_{y,y}(\infty)})}
\end{equation}
where
$\Delta'(0):=\left.\frac{d}{d\lambda}\det(I-D(\lambda))\right|_{\lambda=0+}$.
When $n>1$, $i=1,\ldots,N$ and $y\in S$, the values $b_n(w_i,y)$ are
evaluated by the iterative scheme
\begin{equation}\label{b_n(w_i,w_j)_definition}
b_n(w_i,y)=\sum_{k=1}^N\frac{\alpha_k
f''_k(1)\Delta_{k,i}(0)}{2(n-1)\Delta'(0)}\sum_{r=1}^{n-1}\binom{n}rb_r(w_k,y)b_{n-r}(w_k,y).
\end{equation}
Thus the values $b_n(x,y)$ for $n\geq1$, $x\in S\setminus
W$ and $y\in S$ are determined by the formula
\begin{equation}\label{b_n(x,y)_definition}
b_n(x,y)=\sum\nolimits_{i=1}^N{_{W_i}F_{x,w_i}(\infty)}b_n(w_i,y).
\end{equation}

For $i=1,\ldots,N$ put
\begin{equation}\label{B_1(w_i)_definition}
B_1(w_i)=\sum_{j=1}^N\frac{\alpha_j(f'_j(1)-1)\Delta_{j,i}(0)}{\Delta'(0)}.
\end{equation}
When $n>1$ the values $B_n(w_i)$ are computed iteratively, namely,
\begin{equation}\label{B_n(w_i)_definition}
B_n(w_i)=\sum_{j=1}^N\frac{\alpha_jf''_j(1)\Delta_{j,i}(0)}{2(2n-1)\Delta'(0)}
\sum_{r=1}^{n-1}\binom{n}rB_r(w_j)B_{n-r}(w_j),\; i=1,\ldots,N.
\end{equation}
Bearing on the latter equalities the values $B_n(x)$ for $n\geq1$
and $x\in S\setminus W$ are given as follows
\begin{equation}\label{B_n(x)_definition}
B_n(x)=\sum\nolimits_{i=1}^N{_{W_i}F_{x,w_i}(\infty)}B_n(w_i).
\end{equation}

At last, for subcritical CBP define functions $C_n(x)$, $x\in S$,
$n\in{\bf N}$, arising in the description of the asymptotic behavior for $M_n(t;x)$.
Set
\begin{equation}\label{C_1(w_i)_definition}
C_1(w_i)=1+\sum_{j=1}^N\frac{\alpha_j(f'_j(1)-1)\Delta_{j,i}(0)}{\Delta(0)},\quad
i=1,\ldots,N,
\end{equation}
\begin{equation}\label{C_1(x)_definition}
C_1(x)=1-\sum\nolimits_{i=1}^N{_{W_i}F_{x,w_i}(\infty)}+\sum\nolimits_{i=1}^N{_{W_i}F_{x,w_i}(\infty)C_1(w_i)},\quad
x\in S\setminus W,
\end{equation}
and, for $n>1$, put
\begin{equation}\label{C_n(w_i)_definition}
C_n(w_i)=\sum_{j=1}^N\frac{\Delta_{j,i}(0)}{\Delta(0)}h_{n,j}(C_1(w_j),\ldots,C_{n-1}(w_j)),\quad
i=1,\ldots,N,
\end{equation}
\begin{equation}\label{C_n(x)_definition}
C_n(x)=\sum\nolimits_{i=1}^N{_{W_i}F_{x,w_i}(\infty)}C_n(w_i),\quad
x\in S\setminus W.
\end{equation}

The following theorem is the main result of Section \ref{s:MACBP}
providing the moment analysis of the local and total particles
numbers in CBP. This theorem generalizes Theorems~4.1, 4.2 in
\cite{Y_TPA_11} and Theorem~1 in \cite{DR_13} as well as some
statements of Theorem~2 in \cite{ABY_98}.

\begin{Thm}\label{T:MACBP}
For each $n\in{\bf N}$, given that ${\sf E}\xi^n_i<\infty$ for any
$i=1,\ldots,N$, functions $m_n(t;x,y)$ and $M_n(t;x)$ are bounded on
every finite interval in $[0,\infty)$ for fixed $x,y\in S$.
Moreover, whenever ${\sf E}\xi^n_i<\infty$ for any $i=1,\ldots,N$,
the asymptotic behavior of $m_n(t;x,y)$ and $M_n(t;x)$ as
$t\to\infty$ depends essentially on the class of CBP and is as
follows.
\begin{enumerate}
  \item If $\rho(D)>1$ then $\nu>0$ and
  \begin{eqnarray}
  m_n(t;x,y)&=&a_n(x,y)e^{n\nu t}+o\left(e^{n\nu t}\right),\label{T:MACBP,super_small}\\
  M_n(t;x)&=&A_n(x)e^{n\nu t}+o\left(e^{n\nu t}\right),\quad t\to\infty,\label{T:MACBP,super_big}
  \end{eqnarray}
  where functions $a_n(x,y)$ and $A_n(x)$, $x,y\in S$, are strictly positive.
  \item If $\rho(D)=1$ then
  \begin{eqnarray}
  m_n(t;x,y)&=&b_n(x,y)t^{n-1}+o\left(t^{n-1}\right),\label{T:MACBP,medium_small}\\
  M_n(t;x)&=&B_n(x)t^{2n-1}+o\left(t^{2n-1}\right),\quad t\to\infty,\label{T:MACBP,medium_big}
  \end{eqnarray}
  where function $b_n(x,y)$ is strictly positive for any values $x,y\in S$
  iff ${\int_0^\infty
  u\,d{_{W_j}\overline{F}_{w_i,w_j}(u)}<\infty}$ for each
  $i,j=1,\ldots,N$ and, in addition, in case of recurrent $\eta$ and $n\geq2$ one has
  $\sum\nolimits_{i=1}^N\alpha_i\left|f_i(s)-s\right|>0$ for some $s\in[0,1)$.
  Moreover, $B_n(x)$ is strictly positive for any $x\in S$
  iff the Markov chain $\eta$ is transient and $\int_0^\infty
  u\,d{_{W_j}\overline{F}_{w_i,w_j}(u)}<\infty$ for each
  $i,j=1,\ldots,N$. Otherwise, respectively ${b_n(\cdot,\cdot)\equiv0}$ and
  $B_n(\cdot)\equiv0$.
  \item If $\rho(D)<1$ then
  \begin{eqnarray}
  m_n(t;x,y)&=&o(1),\label{T:MACBP,sub_small}\\
  M_n(t;x)&=&C_n(x)+o(1),\quad t\to\infty,\label{T:MACBP,sub_big}
  \end{eqnarray}
  where function $C_n(x)$ is strictly positive for each $x\in S$ iff
  the Markov chain $\eta$ is transient. Otherwise,
  $C_n(\cdot)\equiv0$.
\end{enumerate}
\end{Thm}
Stress that we do not impose any restrictions on the character of
the Markov chain $\eta$ such as symmetry and homogeneity of
transition probabilities or finite variance of jump sizes etc.
Therefore, in Theorem \ref{T:MACBP} we establish exact asymptotic
relations in item 1 and only partially in items~2 and~3. It is not
possible to give any useful and general information on items~2 and~3
without further knowledge of the underlying motion. For instance,
even for critical symmetric branching random walk on ${\bf Z}^d$
with a single catalyst (see \cite{ABY_98}) there are four different
asymptotic formulae for $m_n(t;x,y)$ and $M_n(t;x)$ depending on
dimension $d=1,2,3$ or $4$ and thereby on the decay rate of
transition probabilities (note that for $d\geq5$ Theorem
\ref{T:MACBP} generalizes the corresponding results in
\cite{ABY_98}). Moreover, for subcritical branching random walk on
${\bf Z}^d$ with a single catalyst (see \cite{B_TrudyMIAN_13}) the
decay orders of $m_n(t;x,y)$ do not coincide for different $d\in{\bf
N}$.

Let us proceed to the proof of Theorem \ref{T:MACBP}.

{\sc Proof.} We will mainly employ the auxiliary Bellman-Harris
process constructed in Section~\ref{s:ABHP}. Firstly we establish
Theorem \ref{T:MACBP} for $x,y\in W$ (Case 1) and afterwards extend
these results to the general case $x,y\in S$ (Case 2).

{\large {\bf Case 1.}} So, assume that $x,y\in W$. The following
proof of Theorem~\ref{T:MACBP} for Case~1 is divided into $9$ steps.
At step 1 we derive a system of renewal equations involving the
factorial moments of the local particles numbers in CBP. Then at
steps 2, 3 and 4 we employ this system to prove the corresponding
assertions (\ref{T:MACBP,super_small}), (\ref{T:MACBP,medium_small})
and (\ref{T:MACBP,sub_small}) of Theorem~\ref{T:MACBP} related to
moment analysis of the local particles numbers. Next at steps 5 and
6 we derive a system of renewal equations involving the factorial
moments of the total particles number in CBP depending on whether
the Markov chain $\eta$ is recurrent or transient. At last, at steps
7, 8 and 9 we exploit these systems to establish statements
(\ref{T:MACBP,super_big}), (\ref{T:MACBP,medium_big}) and
(\ref{T:MACBP,sub_big}) of Theorem~1, respectively, concerning the
moment analysis of the total particles number. Note that at step 6
we introduce an auxiliary decomposable $(L+1)$-dimensional
Bellman-Harris process with a final type of particles.

{\bf Step 1.} At this stage we derive a system of renewal equations
involving the factorial moments of local particles numbers in CBP.
Denote by $m^{BH}_n(t;k,l):={\bf
E}_kZ_l(t)(Z_l(t)-1)\ldots(Z_l(t)-n+1)$, for $t\geq0$, $n\in{\bf N}$
and ${k,l=1,\ldots,L}$, the $n$-th factorial moment of the number of
particles $Z_l(t)$ of the $l$-th type at time $t$ in the
Bellman-Harris process (index $k$ stands for the type of the parent
particle). According to the construction of the auxiliary process
one has $m_n(t;w_i,w_j)=m^{BH}_n(t;i,j)$ for each $n\in{\bf N}$,
$t\geq0$, $i,j=1,\ldots,N$. Thus, given that ${\bf E}\xi^n_k<\infty$
for some $n\in{\bf N}$ and any $k=1,\ldots,N$, the finiteness of
functions $m_n(t;w_i,w_j)$ for each $t\geq0$ and $i,j=1,\ldots,N$,
follows from Theorem~1 in \cite{Sew_74}, Ch.8, Sec.6. Moreover, in
view of this theorem functions $m^{BH}_1(t;k,l)$, $t\geq0$,
$k,l=1,\ldots,L$, satisfy the following system of renewal equations
\begin{equation}\label{m_BH_1(t;k,l)_equation}
m^{BH}_1(t;k,l)=\delta_{k,l}(1-G_k(t))+\sum\limits_{r=1}^Lm_{k,r}\int\nolimits_0^t
m^{BH}_1(t-u;r,l)\,dG_k(u)
\end{equation}
while for $n>1$ this theorem and the generalized Fa\'{a} di Bruno's
formula (see \cite{Good_61}) entail
\begin{eqnarray}\label{m_BH_n(t;k,l)_equation}
& &m^{BH}_n(t;k,l)=\sum_{r=1}^Lm_{k,r}
\int\nolimits_{0}^t{m^{BH}_n(t-u;r,l)\,dG_k(u)}\\
&+&\sum^{\sim}\frac{n!}{\prod\nolimits_{p=1}^{n-1}i^{(1)}_p!\ldots
i^{(L)}_p!}\left.\frac{\partial^{J(1)+\ldots+J(L)}g_k({\bf
s})}{\partial s_1^{J(1)}\ldots\partial s_L^{J(L)}}\right|_{{\bf
s}=(1,\ldots,1)}\nonumber\\
&\times&\int\nolimits_0^t\prod_{p=1}^{n-1}\left(\left(\frac{m^{BH}_p(t-u;1,l)}{p!}\right)^{i^{(1)}_p}
\ldots\left(\frac{m^{BH}_p(t-u;L,l)}{p!}\right)^{i^{(L)}_p}\right)\,dG_k(u).\nonumber
\end{eqnarray}
Here the symbol $\sum\limits^{\sim}$ means the sum taken over all
$L$-dimensional vectors ${{\bf
i}_1=\left(i^{(1)}_1,\ldots,i^{(L)}_1\right)}$, $\ldots$, ${\bf
i}_{n-1}=\left(i^{(1)}_{n-1},\ldots,i^{(L)}_{n-1}\right)$ with
nonnegative integer components satisfying the following equality
$\sum\nolimits_{p=1}^{n-1}p\left(i^{(1)}_p+\ldots+i^{(L)}_p\right)=n$.
For every ${r=1,\ldots,L}$, we also put
${J(r):=\sum\nolimits_{p=1}^{n-1}i^{(r)}_p}$. Substituting the
explicit formulae (\ref{Gg_definition_1}) and
(\ref{Gg_definition_2}) for the offspring generating functions
$g_k(\cdot)$, $k=1,\ldots,L$, of the Bellman-Harris process in
(\ref{m_BH_n(t;k,l)_equation}) we come to more concise relation
\begin{eqnarray}\label{m_BH_n(t;k,l)_equation_new}
& &m^{BH}_n(t;k,l)=\sum_{r=1}^Lm_{k,r}
\int\nolimits_{0}^t{m^{BH}_n(t-u;r,l)\,dG_k(u)}\\
&+&{\bf
I}(k\in\{1,\ldots,N\})\sum\limits^{\approx}\frac{\alpha_kf^{(J)}_k(1)\,n!}{i_1!\ldots
i_{n-1}!}\int\nolimits_0^t
\prod_{p=1}^{n-1}\left(\frac{m^{BH}_p(t-u;k,l)}{p!}\right)^{i_p}\!\!\!dG_k(u)\nonumber
\end{eqnarray}
where ${\bf I}(\cdot)$ is the indicator and the symbol
$\sum\limits^{\approx}$ means the sum taken over all nonnegative
integer $i_1$, $\ldots$, $i_{n-1}$ such that
$\sum\nolimits_{p=1}^{n-1}p\,i_p=n$. We also use the notation
$J:=i_1+\ldots+i_{n-1}$.

{\bf Step 2.} Now we prove the part of assertion
(\ref{T:MACBP,super_small}) related to the case $x,y\in W$. So,
assume that $\rho(D)>1$. Applying Theorem 2.1, item (iii), in
\cite{Crump_70} to the system of renewal equations
(\ref{m_BH_1(t;k,l)_equation}) we get
\begin{equation}\label{m_BH_1_super_sim}
e^{-\nu t}\,m^{BH}_1(t;i,j)\to\frac{(-1)^{i+j}\det(I-H(\nu))_{j,i}}
{(\nu+\beta_j)\left.\frac{d}{d\lambda}\det(I-H(\lambda))\right|_{\lambda=\nu}},\quad
t\to\infty,
\end{equation}
for each $i,j=1,\ldots,N$. Here $\nu>0$ and
$\left.\frac{d}{d\lambda}\det(I-H(\lambda))\right|_{\lambda=\nu}>0$
(see, e.g., \cite{Mode_68_1}). Comparing formulae
(\ref{a_1(w_i,w_j)_definition}) and (\ref{T:MACBP,super_small}) with
(\ref{m_BH_1_super_sim}) we see that one needs to establish
relations between determinants and algebraic adjuncts of matrices
$I-H(\lambda)$ and $I-D(\lambda)$. The following lemma provides
them.
\begin{Lm}\label{L:determinants}
For each $\lambda\geq0$ and $i=1,\ldots,|K_j|$, $j,k=1,\ldots,N$,
one has
$$\det(I-H(\lambda))=\det(I-D(\lambda)),\quad \det(I-H(\lambda))_{j,k}=\det(I-D(\lambda))_{j,k},$$
$$(-1)^{L(j)+i+k}\det(I-H(\lambda))_{L(j)+i,k}=(-1)^{j+k}m_{j,L(j)+i}G^\ast_j(\lambda)\det(I-D(\lambda))_{j,k}.$$
\end{Lm}

{\sc Proof.} We apply transformations to the columns of matrix
$I-H(\lambda)$ which do not change its determinant. Namely, for each
${j,k=1,\ldots,N}$ add the $(L(j)+i)$-th column multiplied by
$G^\ast_{L(j)+i}(\lambda)$ to the $k$-th column if $k\in K_j$, i.e.
$k=k(i,j)$ for some $i=1,\ldots,|K_j|$. After these transformations
we get a block matrix consisting of four blocks. The left upper
block of size $N\times N$ is just $I-D(\lambda)$ whereas the left
lower $(L-N)\times N$ and the right lower $(L-N)\times(L-N)$ blocks
are the zero and the identity matrices, respectively. Employing the
formula for the determinant of a block matrix (see, e.g.,
\cite{Gantmacher_00}, vol.1, Ch.2, Sec.5) we come to the first
assertion of Lemma \ref{L:determinants}. The second one is
established in the same manner. To prove the third assertion we
apply the above transformations to the matrix
$(I-H(\lambda))_{L(j)+i,k}$ and then shift the column number
$L(j)+i-1$ to the column number $N$. The determinants of the matrix
$(I-H(\lambda))_{L(j)+i,k}$ and the transformed one are equal up to
the factor $(-1)^{L(j)+i-N-1}$. Afterwards, we involve the formula
for the determinant of a block matrix once again and use the Laplace
expansion of determinants with respect to the $N$-th column. Lemma
\ref{L:determinants} is proved completely. $\square$

Thus, by virtue of (\ref{m_BH_1_super_sim}) and
Lemma~\ref{L:determinants} we prove (\ref{T:MACBP,super_small}) when
$n=1$ and ${x,y\in W}$. Before focusing on the case $n\!>\!1$ let us
note that in \cite{Mode_68_1} Ch.J.~Mode found the asymptotic
behavior of the mean particles numbers in supercritical multi-type
Bellman-Harris process under rather restrictive conditions that each
c.d.f. $G_k(\cdot)$ has a square-integrable density. This is the
reason why we bear on similar results by K.~Crump established in
\cite{Crump_70} without any additional assumptions and not on the
mentioned results by Ch.J.~Mode. To check relation
(\ref{T:MACBP,super_small}) when $n>1$ and $x,y\in W$ we employ the
induction method in variable $n$. The case $n=1$ is verified above.
Let formula (\ref{T:MACBP,super_small}) be valid for all the moment
orders not exceeding $n-1$. Then the second term at the right-hand
side of (\ref{m_BH_n(t;k,l)_equation_new}) denoted by $V_n(t;k,l)$
has the following asymptotic behavior
$$V_n(t;k,l)\!\sim\!{\bf
I}(k\in\{1,\ldots,N\})e^{n\nu
t}G^\ast_k(n\nu)\sum\limits^{\approx}\frac{\alpha_kf^{(J)}_k(1)\,n!}{i_1!\ldots
i_{n-1}!}\prod_{p=1}^{n-1}\!\!\left(\frac{a_p(w_k,w_l)}{p!}\right)^{i_p}$$
as $t\to\infty$, for each $k=1,\ldots,L$ and $l=1,\ldots,N$. Taking
into account the latter relation and applying Theorem 2.1, item
(iv), in \cite{Crump_70} to the system of renewal equations
(\ref{m_BH_n(t;k,l)_equation_new}) we obtain
$$\frac{m^{BH}_n\!(t;i,j)}{e^{n\nu t}}\!\to\!\!\sum\limits_{k=1}^N\!\!
\frac{(\!-\!1)^{i+k}\beta_k\det(I\!-\!H(n\nu))_{k,i}}
{(n\nu+\beta_k)\det(I\!-\!H(n\nu))}
h_{n,k}(a_1\!(w_k,w_j),\!...,a_{n-\!1}\!(w_k,w_j))$$ as
$t\to\infty$, for each $i,j=1,\ldots,N$. Here we use the alternative
representation (see, e.g., \cite{AB_00}, Theorem~3.3) for function
$h_{n,k}$ defined in (\ref{h_n,k_definition}). The latter asymptotic
relation combined with Lemma~\ref{L:determinants} and formula
(\ref{a_n(w_i,y)_definition}) leads to the desired statement in
(\ref{T:MACBP,super_small}) when $n>1$ and $x,y\in W$. Note that
functions $a_n(w_i,w_j)$, $n\geq1$, $i,j=1,\ldots,N$, are strictly
positive according to, e.g., \cite{Gantmacher_00}, vol.2, Ch.13,
Sec.3.

{\bf Step 3.} Next we prove assertion (\ref{T:MACBP,medium_small})
for $x,y\in W$. Assume that $\rho(D)=1$ or, equivalently, that the
unique solution to $\rho(D(\nu))=1$ is $\nu=0$. Applying Theorem
2.1, item (iii), in \cite{Crump_70} to
(\ref{m_BH_1(t;k,l)_equation}) we come to the following relation
$$m^{BH}_1(t;i,j)\to\frac{(-1)^{i+j}\det(I-H(0))_{j,i}}
{\beta_j\left.\frac{d}{d\lambda}\det(I-H(\lambda))\right|_{\lambda=0+}},\quad
t\to\infty,\quad i,j=1,\ldots,N.$$ This formula along with
(\ref{b_1(w_i,w_j)_definition}) and Lemma~\ref{L:determinants}
entail the required statement of (\ref{T:MACBP,medium_small}) when
$n=1$. For $n>1$ we employ the induction method in variable $n$.
Performing the induction step we get
\begin{eqnarray*}
V_n(t;k,l)&\sim&\sum\limits^{\approx}\frac{\alpha_kf^{(J)}_k(1)\,n!}{i_1!\ldots
i_{n-1}!}\prod_{p=1}^{n-1}\left(\frac{b_p(w_k,w_l)}{p!}\right)^{i_p}\int\nolimits_0^t
(t-u)^{n-J}dG_k(u)\\
&\sim&t^{n-2}\,\frac{\alpha_kf''_k(1)}{2}\sum_{r=1}^{n-1}\binom{n}r
b_r(w_k,w_l)b_{n-r}(w_k,w_l),\quad t\to\infty,
\end{eqnarray*}
for $k,l=1,\ldots,N$, whereas $V_n(t;k,l)=0$, for $k=N+1,\ldots,L$,
$l=1,\ldots,N$. Hence, applying Theorem~2.1, item (v), in
\cite{Crump_70} to (\ref{m_BH_n(t;k,l)_equation_new}) we see that
$$\frac{m^{BH}_n\!(t;i,j)}{t^{n-1}}\!\to\!\!\sum_{k=1}^N
\!\!\frac{(\!-\!1)^{i+k}\!\alpha_kf''_k(1)\!\det(I\!-\!H(0))_{k,i}}{2(n\!-\!1)
\!\!\left.\frac{d}{d\lambda}\det(I\!-\!H(\lambda))\right|_{\lambda=0+}}
\!\sum_{r=1}^{n-1}\!\!\binom{\!n\!}r
b_r\!(w_k,\!w_j)b_{n\!-r}\!(w_k,\!w_j)$$ as $t\to\infty$, for each
$i,j=1,\ldots,N$. The assertion (\ref{T:MACBP,medium_small}) for
$n>1$ now follows from the latter relation, definition
(\ref{b_n(w_i,w_j)_definition}) and Lemma~\ref{L:determinants}. In
view of \cite{Gantmacher_00}, vol.2, Ch.13, Sec.3, the algebraic
adjunct $\Delta_{i,j}(0)$ is strictly positive for any
$i,j=1,\ldots,N$. As was established in \cite{Mode_68_1} one has
$\Delta'(0)\in(0,\infty]$. Moreover, $\Delta'(0)<\infty$ iff
$\int\nolimits_0^\infty
u\,d{_{W_j}\overline{F}_{w_i,w_j}(u)}<\infty$ for each
$i,j=1,\ldots,N$. Therefore, functions $b_n(w_i,w_j)$ are strictly
positive for any $n\in{\bf N}$, $i,j=1,\ldots,N$ if and only if
$\int\nolimits_0^\infty
u\,d{_{W_j}\overline{F}_{w_i,w_j}(u)}<\infty$  and, in addition, in
case of recurrent $\eta$ and $n\geq2$ one has
$\sum\nolimits_{i=1}^N\alpha_i\left|f_i(s)-s\right|>0$ for some
$s\in[0,1)$. Note that the latter condition allows us to separate
the case of ordinary motion of particles without breeding. Such CBP
is critical for recurrent Markov chain $\eta$  and subcritical for
transient $\eta$.

{\bf Step 4.} Relation (\ref{T:MACBP,sub_small}) for the case
$x,y\in W$ is verified by application of Theorem 2.2, item (ii), in
\cite{Mode_68_2} to (\ref{m_BH_1(t;k,l)_equation}) and
(\ref{m_BH_n(t;k,l)_equation_new}). Here we essentially bear on the
fact that $M^n$ converges elementwise to the zero matrix, as
$n\to\infty$, since in subcritical case $\rho(M)<1$ on account of
Lemma~\ref{L:Perron_root}.

{\bf Step 5.} Let us derive a system of renewal equations involving
the factorial moments of the total particles number when $x\in W$
and $\eta$ is recurrent. Denote by $Z(t):=\sum_{j=1}^LZ_j(t)$,
$t\geq0$, the total particles number at time $t$ in the auxiliary
Bellman-Harris process. According to the construction of this
process the distribution laws of $\mu(t)$ and $Z(t)$ are the same
for each $t\geq0$ iff the Markov chain $\eta$ is recurrent. Indeed,
whenever the Markov chain $\eta$ is transient, with positive
probability there are particles at time $t$ in CBP having infinite
life-lengths since they will never hit the set $W$ after time $t$.
These particles are not comprised by our $L$-dimensional
Bellman-Harris process. So, we will now concentrate on the recurrent
case. Then one has $M_n(t;w_i)=M^{BH}_n(t;i)$ for each $t\geq0$,
$n\in{\bf N}$ and $i=1,\ldots,N$ where $M^{BH}_n(t;k):={\bf
E}_kZ(t)(Z(t)-1)\ldots(Z(t)-n+1)$, $t\geq0$, $n\in{\bf N}$ and
${k=1,\ldots,L}$. Thus, given that ${\bf E}\xi^n_l<\infty$ for some
$n\in{\bf N}$ and any ${l=1,\ldots,N}$, Theorem~1 in \cite{Sew_74},
Ch.8, Sec.6, implies the finiteness of functions $M_n(t;w_i)$ for
each $t\geq0$ and $i=1,\ldots,N$. Moreover, in view of this theorem
functions $M^{BH}_1(t;k)$, $t\geq0$, $k=1,\ldots,L$, satisfy the
following system of renewal equations
\begin{equation}\label{M_BH_1(t;k)_equation_recurrent}
M^{BH}_1(t;k)=1-G_k(t)+\sum\limits_{r=1}^Lm_{k,r}\int\nolimits_0^t
M^{BH}_1(t-u;r)\,dG_k(u).
\end{equation}
For $n>1$ this theorem and the generalized Fa\'{a} di Bruno's
formula (see \cite{Good_61}) entail equations in $M^{BH}_n(t;k)$ the
same as those obtained from (\ref{m_BH_n(t;k,l)_equation}) after
replacing $m^{BH}_n(t;k,l)$ by $M^{BH}_n(t;k)$. Therefore, for
$M^{BH}_n(t;k)$, $t\geq0$, $k=1,\ldots,L$, the system of renewal
equations (\ref{m_BH_n(t;k,l)_equation_new}) holds true if function
$m^{BH}_n(t;k,l)$ is replaced by $M^{BH}_n(t;k)$.

{\bf Step 6.} Next consider CBP with the underlying motion of
particles governed by transient Markov chain $\eta$. To cover the
transient case we construct a new auxiliary Bellman-Harris process
with particles of $L+1$ types such that the $(L+1)$-th type is final
(see, e.g., \cite{Sew_74}, Ch.5, Sec.3, and \cite{V_TPA_94}). For
the new Bellman-Harris process, the c.d.f. of the life-length of a
particle of $j$-th type is denoted by $\widehat{G}_j(t)$, $t\geq0$,
whereas its offspring generating function is $\widehat{g}_j({\bf
s},s_{L+1})$, ${{\bf s}\in[0,1]^L}$, $s_{L+1}\in[0,1]$,
$j=1,\ldots,L+1$. Let $\widehat{G}_i(t)=G_i(t)$, $t\geq0$, for
$i=1,\ldots,L$, and $\widehat{g}_j({\bf s},s_{L+1})=g_j({\bf s})$,
${\bf s}\in[0,1]^L$, $s_{L+1}\in[0,1]$, for $j=N+1,\ldots,L$, where
$G_i$ and $g_j$ are defined in (\ref{Gg_definition_1}) and
(\ref{Gg_definition_2}). Then put
$$\widehat{g}_j({\bf s},s_{L+1})=\alpha_j
f_j(s_j)+(1-\alpha_j)\sum\limits_{k=1}^{N}{_{W_k}\overline{F}_{w_j,w_k}(0)s_k}$$
$$+(1-\alpha_j)\sum\limits_{i=1}^{|K_j|}{T_{i,j}(\infty)s_{L(j)+i}}+(1-\alpha_j)\left(1
-\sum\limits_{k=1}^{N}{_{W_k}\overline{F}_{w_j,w_k}(\infty)}\right)s_{L+1},$$
for $j=1,\ldots,N$, ${\bf s}\in[0,1]^L$, $s_{L+1}\in[0,1]$, i.e.
definitions of functions $g_j$ and $\widehat{g}_j$ differ by the
last term only. Since the $(L+1)$-th type is final, that is every
particle of the $(L+1)$-th type has infinite life-length and does not
produce any offsprings, for the sake of definiteness we may set
${\widehat{G}_{L+1}(t)=0}$, $t\geq0$, and $\widehat{g}_{L+1}({\bf
s},s_{L+1})=s_{L+1}$, ${\bf s}\in[0,1]^L$, ${s_{L+1}\in[0,1]}$. The
mean matrix $\widehat{M}=(\widehat{m}_{k,l})_{k,l=1}^{L+1}$ of the
new process has entries ${\widehat{m}_{k,l}=m_{k,l}}$, for
$k,l=1,\ldots,L$, $\widehat{m}_{i,L+1}=(1-\alpha_i)\left(1
-\sum\nolimits_{k=1}^{N}{_{W_k}\overline{F}_{w_i,w_k}(\infty)}\right)$,
$i=1,\ldots,N$, $\widehat{m}_{L+1,L+1}=1$ and $\widehat{m}_{k,l}=0$
for the remaining pairs of $k$ and $l$. Let $\widehat{Z}_j(t)$,
$t\geq0$, ${j=1,\ldots,L+1}$, be the number of particles of type $j$
at time~$t$ in the new Bellman-Harris process. Evidently, the
distributions of vectors
$(\widehat{Z}_1(t),\ldots,\widehat{Z}_L(t))$ and
$(Z_1(t),\ldots,Z_L(t))$ coincide for each $t\geq0$ whenever the
parent particles of both Bellman-Harris processes have the same
type. Moreover, the new process takes into account even the
particles from the ${(L+1)}$-th group of particles in CBP (see
Section~\ref{s:ABHP}) which are not comprised by the $L$-dimensional
Bellman-Harris process. Thus, for each $n\in{\bf N}$,
${i=1,\ldots,N}$ and $t\geq0$, one has
$M_n(t;w_i)=\widehat{M}^{BH}_n(t;i)$ where we write
${\widehat{M}^{BH}_n(t;k):={\bf
E}_k\widehat{Z}(t)(\widehat{Z}(t)-1)\ldots(\widehat{Z}(t)-n+1)}$,
$k=1,\ldots,L+1$, and set
${\widehat{Z}(t):=\sum_{j=1}^{L+1}\widehat{Z}_j(t)}$. Then applying
Theorem~1 in \cite{Sew_74}, Ch.8, Sec.6, we ascertain that whenever
${\bf E}\xi^n_l<\infty$ for some $n\in{\bf N}$ and any
$l=1,\ldots,N$, functions $M_n(t;w_i)$ are finite for each $t\geq0$
and $i=1,\ldots,N$. This theorem also implies that functions
$\widehat{M}^{BH}_1(t;k)$, $t\geq0$, $k=1,\ldots,L+1$, satisfy the
following system of $L+1$ renewal equations
$$\widehat{M}^{BH}_1(t;k)=1-\widehat{G}_k(t)+\sum\limits_{r=1}^{L+1}\widehat{m}_{k,r}\int\nolimits_0^t
\widehat{M}^{BH}_1(t-u;r)\,d\widehat{G}_k(u).$$ The $(L+1)$-th
equation is just $\widehat{M}^{BH}_1(t;L+1)=1$, $t\geq0$.
Substituting that value $1$ and expressions for $\widehat{m}_{k,r}$,
$r=1,\ldots,L+1$,
into the first $L$ equations we get a new system of
$L$ renewal equations
\begin{eqnarray}\label{M_BH_1(t;k)_equation_transient}
&
\!\!&\!\!\widehat{M}^{BH}_1(t;k)=\sum\limits_{r=1}^{L}m_{k,r}\int\nolimits_0^t
\widehat{M}^{BH}_1(t-u;r)\,dG_k(u)\\
&+\!\!&\!\!1-G_k(t)+{\bf
I}\left(k\in\{1,\ldots,N\}\right)(1-\alpha_k)\left(1
-\sum\limits_{j=1}^{N}{_{W_j}\overline{F}_{w_k,w_j}(\infty)}\right)G_k(t)\nonumber
\end{eqnarray}
for $t\geq0$ and $k=1,\ldots,L$. It is not difficult to verify that
in view of Theorem~1 in \cite{Sew_74}, Ch.8, Sec.6, the system of
$L+1$ equations resulting from
(\ref{m_BH_n(t;k,l)_equation}) upon replacement of $L$ by $L+1$ and
of $m^{BH}_n(t;k,l)$ by $\widehat{M}^{BH}_n(t;k)$ holds true for
$n>1$, $t\geq0$ and $k=1,\ldots,L+1$. The $(L+1)$-th equation is just
$\widehat{M}^{BH}_n(t;L+1)=0$. Substituting that in the first $L$
equations and simplifying them similar to derivation of
(\ref{m_BH_n(t;k,l)_equation_new}) we get the system of $L$ renewal
equations in functions $\widehat{M}^{BH}_n(t;k)$, for $n>1$,
$t\geq0$ and $k=1,\ldots,L$, analogous to
(\ref{m_BH_n(t;k,l)_equation_new}) but with $m^{BH}_n(t;k,l)$
replaced by $\widehat{M}^{BH}_n(t;k)$.

{\bf Step 7.} Now we prove assertion (\ref{T:MACBP,super_big}) for
the case $x,y\in W$. So, assume that $\rho(D)>1$. At first consider
the recurrent case. Applying Theorem 2.1, item (iii), in
\cite{Crump_70} to (\ref{M_BH_1(t;k)_equation_recurrent}) we see
that
\begin{eqnarray}\label{M_BH_1_super_sim_reccurent}
& &\frac{M^{BH}_1(t;i)}{e^{\nu
t}}\!\to\!\!\sum_{j=1}^N\!\frac{\Delta_{j,i}(\nu)}{\nu\Delta'(\nu)}\!\!
\left(\!1-\alpha_jG^\ast_j(\nu)-(1-\alpha_j)G^\ast_j(\nu)\!\!\sum_{k=1}^N
\!{_{W_k}\overline{F}^\ast_{w_j,w_k}(\nu)}\!\!\right)\\
&=&\sum_{k=1}^N\sum_{j=1}^N\frac{\Delta_{j,i}(\nu)}{\nu\Delta'(\nu)}
\left(\delta_{j,k}-\delta_{j,k}\alpha_jf'_j(1)G^\ast_j(\nu)-(1-\alpha_j)G^\ast_j(\nu)
{_{W_k}\overline{F}^\ast_{w_j,w_k}(\nu)}\right)\nonumber\\
&+&\sum_{j=1}^N\frac{\alpha_j(f'_j(1)-1)G^\ast_j(\nu)\Delta_{j,i}(\nu)}{\nu\Delta'(\nu)}=
\sum_{j=1}^N\frac{\alpha_j\beta_j(f'_j(1)-1)\Delta_{j,i}(\nu)}{\nu(\nu+\beta_j)\Delta'(\nu)}\nonumber
\end{eqnarray}
as $t\to\infty$ and $i=1,\ldots,N$. Here we rely on
Lemma~\ref{L:determinants} and formulae
(\ref{Gg_definition_1})-(\ref{matrix_M=}) as well as on the Laplace
expansion of determinants and obvious equality
${\sum_{k=1}^N{_{W_k}\overline{F}_{w_j,w_k}(\infty)}=1}$ valid in
the recurrent case. We also take into account that $\Delta(\nu)=0$
since $1$ is an eigenvalue of matrix $D(\nu)$. Relation
(\ref{M_BH_1_super_sim_reccurent}) along with
(\ref{A_1(w_i)_definition}) lead to the desired statement
(\ref{T:MACBP,super_big}) when $n=1$ and the Markov chain $\eta$ is
recurrent. The transient case is treated in the same way except for
two differences. The first one consists in employment of equations
system (\ref{M_BH_1(t;k)_equation_transient}) instead of
(\ref{M_BH_1(t;k)_equation_recurrent}). The second one is that for
transient Markov chain $\eta$ the strict inequality
${\sum_{k=1}^N{_{W_k}\overline{F}_{w_j,w_k}(\infty)}<1}$ holds true
at least for some $j\in\{1,\ldots,N\}$. However, when deriving
relation for $\widehat{M}^{BH}_1(t;i)$ similar to
(\ref{M_BH_1_super_sim_reccurent}), the gap between the latter sum
and $1$ is compensated by means of the additional summand in
(\ref{M_BH_1(t;k)_equation_transient}). Thus, when $n=1$ assertion
(\ref{T:MACBP,super_big}) is established and the asymptotic behavior
of $M_1(t;w_i)$, $i=1,\ldots,N$, in supercritical CBP does not
depend on whether the Markov chain $\eta$ is recurrent or transient.
Since for $n>1$ the systems of equations  in $M^{BH}_n(t;k)$ and
$\widehat{M}^{BH}_n(t;k)$ respectively are of the same type as
(\ref{m_BH_n(t;k,l)_equation_new}), the asymptotic behavior of
$M_n(t;w_i)$ in supercritical CBP is also independent of recurrence
or transience of $\eta$. Moreover, when $n>1$ derivation of
(\ref{T:MACBP,super_big}) for the total particles number almost
literally repeats that of (\ref{T:MACBP,sub_small}) for the local
particles numbers and thus is omitted. Note that function $A_1(w_i)$
(and, consequently, $A_n(w_i)$ for $n>1$) is strictly positive for
each $i=1,\ldots,N$ in view of positivity of the terms at the
right-hand side of (\ref{M_BH_1_super_sim_reccurent}).

{\bf Step 8.} Next proceed to the proof of statement
(\ref{T:MACBP,medium_big}) for the case $x,y\in W$. So, consider the
case $\rho(D)=1$. Applying Theorem 2.1, item (v), in \cite{Crump_70}
to (\ref{M_BH_1(t;k)_equation_transient}) we get
\begin{eqnarray}
&&t^{-1}\widehat{M}^{BH}_1(t;i)\to\sum_{j=1}^N\frac{(1-\alpha_j)\Delta_{j,i}(0)}{\Delta'(0)}
\left(1-\sum_{k=1}^N{_{W_k}\overline{F}_{w_j,w_k}(\infty)}\right)\label{M_BH_1_medium_sim_transient}\\
&=&\sum_{k=1}^N\sum_{j=1}^N\frac{\Delta_{j,i}(0)}{\Delta'(0)}
\left(\delta_{j,k}-\delta_{j,k}\,\alpha_jf'_j(1)-(1-\alpha_j){_{W_k}\overline{F}_{w_j,w_k}(\infty)}\right)
\nonumber\\
&+&\sum_{j=1}^N\frac{\alpha_j(f'_j(1)-1)\Delta_{j,i}(0)}{\Delta'(0)}
=\sum_{j=1}^N\frac{\alpha_j(f'_j(1)-1)\Delta_{j,i}(0)}{\Delta'(0)},\quad
t\to\infty,\nonumber
\end{eqnarray}
for $i=1,\ldots,N$. Here we use Lemma~\ref{L:determinants} as well
as the Laplace expansion of determinants and equality $\Delta(0)=0$
valid by virtue of assumption ${\rho(D)=1}$. Relation
(\ref{M_BH_1_medium_sim_transient}) combined with
(\ref{B_1(w_i)_definition}) implies the desired assertion in
(\ref{T:MACBP,medium_big}) when $n=1$ and Markov chain $\eta$ is
transient. When Markov chain $\eta$ is recurrent, the corresponding
assertion follows from Theorem 2.1, item (i), in \cite{Crump_70} and
observation that the non-integral term at the right-hand side in
(\ref{M_BH_1(t;k)_equation_recurrent}) tends to $0$, as
$t\to\infty$, whereas that in (\ref{M_BH_1(t;k)_equation_transient})
converges to a positive limit at least for some
$k\in\{1,\ldots,N\}$. For $n>1$ derivation of
(\ref{T:MACBP,medium_big}) for the total particles number again
reproduces that of (\ref{T:MACBP,medium_small}) for the local
particles numbers and so is omitted. Note that in view of
(\ref{M_BH_1_medium_sim_transient}) function $B_1(w_i)$ (and,
consequently, $B_n(w_i)$ for $n>1$) is strictly positive for each
$i=1,\ldots,N$ iff $\Delta'(0)<\infty$ and
$\sum_{k=1}^N{_{W_k}\overline{F}_{w_j,w_k}(\infty)}<1$ at least for
some $j\in\{1,\ldots,N\}$.

{\bf Step 9.} At last, let us establish (\ref{T:MACBP,sub_big}) for
the case $x,y\in W$. Thus, assume that $\rho(D)<1$. Applying Theorem
2.2, item (ii), in \cite{Mode_68_2} to
(\ref{M_BH_1(t;k)_equation_transient}), we see that
\begin{eqnarray}
&&\widehat{M}^{BH}_1(t;i)\to\sum_{j=1}^N\frac{(1-\alpha_j)\Delta_{j,i}(0)}{\Delta(0)}
\left(1-\sum_{k=1}^N{_{W_k}\overline{F}_{w_j,w_k}(\infty)}\right)\label{M_BH_1_sub_sim_transient}\\
&=&\sum_{k=1}^N\sum_{j=1}^N\frac{\Delta_{j,i}(0)}{\Delta(0)}
\left(\delta_{j,k}-\delta_{j,k}\,\alpha_jf'_j(1)-(1-\alpha_j){_{W_k}\overline{F}_{w_j,w_k}(\infty)}\right)
\nonumber\\
&+&\sum_{j=1}^N\frac{\alpha_j(f'_j(1)-1)\Delta_{j,i}(0)}{\Delta(0)}
=1+\sum_{j=1}^N\frac{\alpha_j(f'_j(1)-1)\Delta_{j,i}(0)}{\Delta(0)},\quad
t\to\infty,\nonumber
\end{eqnarray}
for each $i=1,\ldots,N$. Here we involve Lemma~\ref{L:determinants}
as well as the Laplace expansion of determinants and inequality
$\Delta(0)>0$ valid on account of assumption $\rho(D)<1$. Formulae
(\ref{C_1(w_i)_definition}) and (\ref{M_BH_1_sub_sim_transient})
entail the desired statement (\ref{T:MACBP,sub_big}) when $n=1$ and
Markov chain $\eta$ is transient. When Markov chain $\eta$ is
recurrent, the corresponding relation ensues from
(\ref{M_BH_1(t;k)_equation_recurrent}) and Theorem~2.2, item (ii),
in \cite{Mode_68_2}. When $n>1$, this theorem also implies
(\ref{T:MACBP,sub_big}) due to equations system
(\ref{m_BH_n(t;k,l)_equation_new}). Observe that in view of
(\ref{M_BH_1_sub_sim_transient}) function $C_1(w_i)$ (and,
consequently, $C_n(w_i)$ for $n>1$) is strictly positive for each
$i=1,\ldots,N$ iff
$\sum_{k=1}^N{_{W_k}\overline{F}_{w_j,w_k}(\infty)}<1$ for at least
some $j\in\{1,\dots,N\}$, i.e. Markov chain $\eta$ is transient.

Thus, Theorem \ref{T:MACBP} is proved completely for the case
$x,y\in W$.

{\large{\bf Case 2.}} Now we assume that either $x\in S\setminus W$
or $y\in S\setminus W$ or ${x,y\in S\setminus W}$. The main idea of
the rest of the proof is as follows. If $x\in S\setminus W$ and
$y\in W$ we supplement the set of catalysts $W$ with $x$. Vice
versa, if $x\in W$ and $y\in S\setminus W$ let us enlarge the set of
catalysts $W$ by $y$. If ${x\in S\setminus W}$ and $y=x$ then we
will add $x$ to the set of catalysts $W$. If both $x$ and $y$ are
from $S\setminus W$ and, moreover, $x\neq y$, we supplement the set
of catalysts $W$ with both states $x$ and $y$. Afterwards we may
employ the already established results for CBP with $N+1$ or $N+2$
catalysts. So, set $w_{N+1}=x$, ${W(x):=W\cup\{x\}}$ and
${W_i(x):=W(x)\setminus\{w_i\}}$, ${i=1,\ldots,N+1}$. Let
${D(x;\lambda)=(d_{i,j}(x;\lambda))_{i,j=1}^{N+1}}$ be matrix with
$d_{i,j}(x;\lambda):=\delta_{i,j}\,\alpha_i
f_i'(1)G^\ast_i(\lambda)+(1-\alpha_i)G^\ast_i(\lambda){_{W_j(x)}\overline{F}^\ast_{w_i,w_j}(\lambda)}$,
$\lambda\geq0$. Here $\alpha_{N+1}=0$, ${f'_{N+1}(1)=0}$ and
$G_{N+1}(t):=1-e^{q(x,x)t}$, ${t\geq0}$. Similarly put $w_{N+2}=y$,
$W(x,y):=W(x)\cup\{y\}$ and ${W_i(x,y):=W(x,y)\setminus\{w_i\}}$,
${i=1,\ldots,N+2}$. Consider matrix
$D(x,y;\lambda)=(d_{i,j}(x,y;\lambda))_{i,j=1}^{N+2}$, for
$\lambda\geq0$, such that
$d_{i,j}(x,y;\lambda):=\delta_{i,j}\,\alpha_i
f_i'(1)G^\ast_i(\lambda)+(1-\alpha_i)G^\ast_i(\lambda){_{W_j(x,y)}\overline{F}^\ast_{w_i,w_j}(\lambda)}$.
Here $\alpha_{N+2}=0$, ${f'_{N+2}(1)=0}$ and
$G_{N+2}(t):=1-e^{q(y,y)t}$, $t\geq0$. Denote by $\nu(x)$ and
$\nu(x,y)$ the unique solutions to equations
$\rho\left(D(x;\nu(x))\right)=1$ and
$\rho\left(D(x,y;\nu(x,y))\right)=1$, respectively. The following
Lemma~\ref{L:3D} implies that $\nu=\nu(x)=\nu(x,y)$.
\begin{Lm}\label{L:3D}
For each $\lambda\geq0$, the Perron roots $\rho(D(\lambda))$,
$\rho(D(x;\lambda))$ and $\rho(D(x,y;\lambda))$ are simultaneously
greater than $1$, equal to $1$ or less than $1$.
\end{Lm}

{\sc Proof.} The arguments of Lemma \ref{L:3D} resemble those of
Lemma~\ref{L:Perron_root}. It is enough to establish the desired
statement for matrices $D(\lambda)$ and $D(x;\lambda)$ only since
automatically the similar assertion is true for $D(x;\lambda)$ and
$D(x,y;\lambda)$. According to the Perron-Frobenius theorem matrix
$D(x;\lambda)$ has a strictly positive left eigenvector ${\bf
v}(x)=\left(v_1(x),\ldots,v_{N+1}(x)\right)^\top$ corresponding to
the eigenvalue $\rho(D(x;\lambda))$, i.e.
${\left(D(x;\lambda)^\top-\rho(D(x;\lambda))I\right){\bf v}(x)={\bf
0}}$. We will apply the  equivalent transformations to the obtained
system of equations, namely, for each ${i=1,\ldots,N}$ add the
$(N+1)$-th row multiplied by ${_{W_i}F^\ast_{w_{N+1},w_i}(\lambda)}$
to the $i$-th row. Focusing on the first $N$ equations only, we
deduce that
\begin{eqnarray}\label{2D}
&
&\left(v_1,\ldots,v_N\right)\left(D(\lambda)-\rho(D(x;\lambda))I\right)\\
&=&v_{N+1}\left(\rho(D(x;\lambda))-1\right)
\left({_{W_1}F^\ast_{w_{N+1},w_1}(\lambda)},\ldots,{_{W_N}F^\ast_{w_{N+1},w_N}(\lambda)}\right).\nonumber
\end{eqnarray}
Here we use the following identity
\begin{equation}\label{taboo_identity}
{_{W_i}\overline{F}^\ast_{w_j,w_i}(\lambda)}={_{W_i(x)}\overline{F}^\ast_{w_j,w_i}(\lambda)}+{_{W_{N+1}(x)}\overline{F}^\ast_{w_j,w_{N+1}}(\lambda)}
{_{W_i}F^\ast_{w_{N+1},w_i}(\lambda)}
\end{equation}
where $i=1,\ldots,N$ and $j=1,\ldots,N+1$. This equality is true,
e.g., in view of Theorem~8 in \cite{Chung_60}, Ch.2, Sec.11. It
follows from (\ref{2D}) that
\begin{eqnarray}
\quad\left(v_1,\ldots,v_N\right)D(\lambda)
\geq\rho(D(x;\lambda))\left(v_1,\ldots,v_N\right)&
\mbox{for}&\rho(D(x;\lambda))\geq1,\label{2D>=}\\
\quad\left(v_1,\ldots,v_N\right)D(\lambda)
\leq\rho(D(x;\lambda))\left(v_1,\ldots,v_N\right)&
\mbox{for}&\rho(D(x;\lambda))\leq1.\label{2D<=}
\end{eqnarray}
Examining the proof of the Perron-Frobenius theorem we conclude that
due to (\ref{2D>=}) one has $\rho(D(\lambda))\geq\rho(D(x;\lambda))$
for $\rho(D(x;\lambda))\geq1$. On account of Theorem 1.6 in
\cite{Seneta_06}, Ch.1, Sec.4, relation (\ref{2D<=}) yields that
${\rho(D(\lambda))\leq\rho(D(x;\lambda))}$ for
$\rho(D(x;\lambda))\leq1$. These estimates entail the desired
assertion of Lemma~\ref{L:3D}. $\square$

After application of the established part of Theorem~\ref{T:MACBP}
to the CBP with the catalysts set $W(x)$ or $W(x,y)$ we observe that
to complete the proof of Theorem~\ref{T:MACBP} one could bring the
expressions like $a_n(w_{N+1},w_{N+2})$ into the form
(\ref{a_n(x,y)_definition}). This can be easily realized by using
the following result.
\begin{Lm}\label{L:x,y_determinants}
For $x,y\in S\setminus W$, $x\neq y$, $i,j=1,\ldots,N$ and
$\lambda\geq0$, one has
\begin{equation}\label{L:x,y_determinants_1}
\det\left(I-D(x;\lambda)\right)=\Delta(\lambda)\left(1-{_WF^\ast_{x,x}(\lambda)}\right),
\end{equation}
\begin{equation}\label{L:x,y_determinants_2}
\det\left(I-D(x,y;\lambda)\right)=
\Delta(\lambda)\left(1-{_WF^\ast_{x,x}(\lambda)}\right)\left(1-{_{W(x)}F^\ast_{y,y}(\lambda)}\right),
\end{equation}
\begin{equation}\label{L:x,y_determinants_3}
\left.\frac{d}{d\lambda}\det\left(I-D(x;\lambda)\right)\right|_{\lambda=\nu}=
\Delta'(\nu)\left(1-{_WF^\ast_{x,x}(\nu)}\right),
\end{equation}
\begin{equation}\label{L:x,y_determinants_4}
\left.\frac{d}{d\lambda}\det\left(I-D(x,y;\lambda)\right)\right|_{\lambda=\nu}
\!=\!\Delta'(\nu)\left(1\!-\!{_WF^\ast_{x,x}(\nu)}\right)\!\left(1\!-\!{_{W(x)}F^\ast_{y,y}(\nu)}\right)\!,
\end{equation}
\begin{equation}\label{L:x,y_determinants_5}
(-1)^{i+j}\det\left(I-D(x;\lambda)\right)_{j,i}=\Delta_{j,i}(\lambda)\left(1-{_WF^\ast_{x,x}(\lambda)}\right),
\end{equation}
\begin{equation}\label{L:x,y_determinants_6}
\frac{(-1)^{j+N+1}\det\left(I-D(x;\lambda)\right)_{j,N+1}}{1-{_WF^\ast_{x,x}(\lambda)}}
=\sum_{i=1}^{N}{_{W_i}F^\ast_{x,w_i}(\lambda)}\Delta_{j,i}(\lambda),
\end{equation}
\begin{equation}\label{L:x,y_determinants_7}
(-1)^{i+N+1}\det\left(I\!-\!D(x;\lambda)\right)_{N+1,i}\!=\!
\sum_{j=1}^{N}\Delta_{j,i}(\lambda)(1\!-\!\alpha_j)G^\ast_j(\lambda){_W\overline{F}^\ast_{w_j,x}(\lambda)},
\end{equation}
\begin{eqnarray}
&
&
\det\!\left(I\!-\!D(x;\lambda)\right)_{\!N\!+\!1,N\!+\!1}\label{L:x,y_determinants_8}\\
&
&=\Delta(\lambda)+\!\!\sum_{i,j=1}^{N}\!\!{_{W_i}F^\ast_{x,w_i}(\lambda)}
\Delta_{j,i}(\lambda)(1\!-\!\alpha_j)G^\ast_j(\lambda){_W\overline{F}^\ast_{w_j,x}(\lambda)},\nonumber
\end{eqnarray}
\begin{eqnarray}
&
&\frac{-\det(I-D(x,y;\lambda))_{N+2,N+1}\left(1-{_WF^\ast_{y,y}(\lambda)}\right)}
{\left(1-{_WF^\ast_{x,x}(\lambda)}\right)\left(1-{_{W(x)}F^\ast_{y,y}(\lambda)}\right)}
\label{L:x,y_determinants_9}\\
&
&=\Delta(\lambda){_WF^\ast_{x,y}(\lambda)}+\sum_{i,j=1}^{N}{_{W_i}F^\ast_{x,w_i}(\lambda)}
\Delta_{j,i}(\lambda)(1-\alpha_j)G^\ast_j(\lambda){_W\overline{F}^\ast_{w_j,y}(\lambda)}\nonumber
\end{eqnarray}
where $\nu\geq0$ is the unique solution to $\rho(D(\nu))=1$ existing
whenever ${\rho(D)\geq1}$. If $\nu=0$ then
$\left.\frac{d}{d\lambda}\det\left(I-D(x;\lambda)\right)\right|_{\lambda=\nu}$
means the right derivative at $0$.
\end{Lm}

{\sc Proof.} Apply some transformations to the columns of matrix
$I-D(x;\lambda)$ which do not change its determinant. Namely, for
each ${i=1,\ldots,N}$ add the $(N+1)$-th column multiplied by
${_{W_i}F^\ast_{w_{N+1},w_i}(\lambda)}$ to the $i$-th column.
Afterwards using (\ref{taboo_identity}) and the Laplace expansion of
the determinants with respect to the ${(N+1)}$-th row we come to
(\ref{L:x,y_determinants_1}). Relation (\ref{L:x,y_determinants_2})
can be considered as ensuing from (\ref{L:x,y_determinants_1}).
Formulae (\ref{L:x,y_determinants_1}) and
(\ref{L:x,y_determinants_2}) imply (\ref{L:x,y_determinants_3}) and
(\ref{L:x,y_determinants_4}) since $\Delta(\lambda)>0$ for
$\lambda>\nu$ and $\Delta(\nu)=0$. Equality
(\ref{L:x,y_determinants_5}) is verified in the same way as
(\ref{L:x,y_determinants_1}).

Let us check the validity of (\ref{L:x,y_determinants_6}).
Obviously, the following identity is true
\begin{equation}\label{taboo_identity_new}
\frac{_{W_k(x)}F^\ast_{x,w_k}(\lambda)}{1-{_WF^\ast_{x,x}(\lambda)}}={_{W_k}F^\ast_{x,w_k}(\lambda)},
\quad\lambda\geq0,\quad k=1,\ldots,N.
\end{equation}
Apply some transformations to the rows of matrix
$\left(I-D(x;\lambda)\right)_{j,N+1}$ which do not change its
determinant. Namely, for each ${i=1,\ldots,N-1}$, add the $N\mbox{-th}$ row
multiplied by
$(1-\alpha_i)G^\ast_i(\lambda){_W\overline{F}^\ast_{w_i,x}(\lambda)}/\left(1-{_WF^\ast_{x,x}(\lambda)}\right)$
to the $i\mbox{-th}$ row. Afterwards employing (\ref{taboo_identity}),
(\ref{taboo_identity_new}) and the Laplace expansion of the
determinants with respect to the $N$-th row we derive
(\ref{L:x,y_determinants_6}).

Implementing the same transformations as for proving
(\ref{L:x,y_determinants_1}) and using the Laplace expansion of the
determinants with respect to the $N$-th column lead to identity
(\ref{L:x,y_determinants_7}). Formula (\ref{L:x,y_determinants_8})
is verified by involving relation (\ref{taboo_identity}), the
linearity property of a determinant with respect to rows and the
Laplace expansion of the determinants.

Relation (\ref{L:x,y_determinants_9}) is established by the
combination of the previous transformations which do not change the
determinant of $(I-D(x,y;\lambda))_{N+2,N+1}$. Firstly, similar to
proving (\ref{L:x,y_determinants_1}), for each ${i=1,\ldots,N}$, add
the $(N+1)$-th column multiplied by
${_{W_i(x)}F^\ast_{w_{N+2},w_i}(\lambda)}$ to the $i$-th column and
take into account a counterpart of (\ref{taboo_identity}). Secondly,
by analogy with verifying (\ref{L:x,y_determinants_6}), for every
${i=1,\ldots,N}$, add to the $i$-th row the $(N+1)$-th row multiplied
by
$(1-\alpha_i)G^\ast_i(\lambda){_W\overline{F}^\ast_{w_i,x}(\lambda)}/\left(1-{_WF^\ast_{x,x}(\lambda)}\right)$.
Thirdly, apply the Laplace expansion of the determinants with
respect to the $(N+1)$-th column and the $(N+1)$-th row. At last,
employ formulae (\ref{taboo_identity}) and
(\ref{taboo_identity_new}) as well as the following equalities
$${_{W(x)}\overline{F}^\ast_{w_i,y}(\lambda)}
+\frac{{_W\overline{F}^\ast_{w_i,x}(\lambda)}{_{W(x)}
F^\ast_{x,y}(\lambda)}}{1-{_WF^\ast_{x,x}(\lambda)}}={_W\overline{F}^\ast_{w_i,y}(\lambda)}
\frac{1-{_{W(x)}F^\ast_{y,y}(\lambda)}}{1-{_WF^\ast_{y,y}(\lambda)}},$$
$$\frac{{_{W(x)}}F^\ast_{x,y}(\lambda)\left(1-{_WF^\ast_{y,y}(\lambda)}\right)}
{\left(1-{_WF^\ast_{x,x}(\lambda)}\right)\left(1-{_{W(x)}F^\ast_{y,y}(\lambda)}\right)}={_WF^\ast_{x,y}(\lambda)},
\quad \lambda\geq0,\quad i=1,\ldots,N,$$ valid according to the
proof of Theorem~2 in \cite{B_SPL_13}. In result we come to
(\ref{L:x,y_determinants_9}). This completes the proof of
Lemma~\ref{L:x,y_determinants}. $\square$

Thus, using Lemma~\ref{L:x,y_determinants} we derive relations
(\ref{a_1(w_i,y)_definition}), (\ref{a_n(x,y)_definition}),
(\ref{A_n(x)_definition}), (\ref{b_n(x,y)_definition}),
(\ref{B_n(x)_definition}), (\ref{C_1(x)_definition}),
(\ref{C_n(x)_definition}) from (\ref{a_1(w_i,w_j)_definition}),
(\ref{a_n(w_i,y)_definition}), (\ref{A_1(w_i)_definition}),
(\ref{A_n(w_i)_definition}), (\ref{b_1(w_i,w_j)_definition}),
(\ref{b_n(w_i,w_j)_definition}), (\ref{B_1(w_i)_definition}),
(\ref{B_n(w_i)_definition}), (\ref{C_1(w_i)_definition}) and
(\ref{C_n(w_i)_definition}). Theorem~\ref{T:MACBP} is proved
completely. $\square$

\section{Applications}\label{s:Applications}

In section \ref{s:MACBP} by virtue of Theorem~\ref{T:MACBP} we
established that the Malthusian parameter $\nu$ plays a crucial role
in the asymptotic behavior of both total and local particles numbers
in supercritical CBP. Before considering some particular examples in
the present section let us discuss the alternative way of evaluating
$\nu$. For this purpose introduce matrices
$\widetilde{D}(\lambda)=\left(\widetilde{d}_{i,j}(\lambda)\right)_{i,j=1}^N$
and
$\widehat{D}(\lambda)=\left(\widehat{d}_{i,j}(\lambda)\right)_{i,j=1}^N$
with the corresponding entries
\begin{eqnarray*}
\widetilde{d}_{i,j}(\lambda)&:=&\left(\alpha_if_i'(1)G^\ast_i(\lambda)-1
+(1-\alpha_i)G^\ast_i(\lambda)\frac{\lambda-q(w_i,w_i)}{-q(w_i,w_i)}\right)F^\ast_{w_i,w_j}(\lambda)\\
&+&
\delta_{i,j}\left(\alpha_if_i'(1)G^\ast_i(\lambda)-1\right)\left(1-F^\ast_{w_i,w_i}(\lambda)\right),\quad
\lambda\geq0,\\
\widehat{d}_{i,j}(\lambda)&:=&\left(\alpha_if_i'(1)G^\ast_i(\lambda)-1
+(1-\alpha_i)G^\ast_i(\lambda)\frac{\lambda-q(w_i,w_i)}{-q(w_i,w_i)}\right)G_\lambda(w_i,w_j)\\
&+&\delta_{i,j}\frac{1-\alpha_i}{q(w_i,w_i)}G^\ast_i(\lambda),\quad\lambda>0.
\end{eqnarray*}
Here $G_\lambda(x,y):=\int\nolimits_{0}^\infty{e^{-\lambda
t}p(t;x,y)\,dt}$, $\lambda>0$, is the Laplace transform of the
transition probability $p(t;x,y)$, $x,y\in S$, $t\geq0$, of the
Markov chain $\eta$. Note that function
$G_0(x,y)=\lim\nolimits_{\lambda\to0+}G_{\lambda}(x,y)$ is called
\emph{Green's function} and is finite iff the Markov chain $\eta$ is
transient (see, e.g., Theorem~4 and Corollary~2 in \cite{Chung_60},
Ch.2, Sec.10). Recall also that according to Theorem~3 and relation
(4) in \cite{Chung_60}, Ch.2, Sec.12, one has
\begin{equation}\label{from_Chung}
F^\ast_{w_i,w_j}(\lambda)=\frac{G_{\lambda}(w_i,w_j)}{G_{\lambda}(w_j,w_j)},\quad
F^\ast_{w_i,w_i}(\lambda)=1-\frac{1}{(\lambda-q(w_i,w_i))G_{\lambda}(w_i,w_i)},
\end{equation}
for any $i,j=1,\ldots,N$, $i\neq j$ and $\lambda>0$.

It follows from Lemma~\ref{L:Perron_root} and definition of the
Malthusian parameter in section \ref{s:MACBP} that $\nu$ can be
found as the maximal $\lambda$ satisfying the relation
$\det(D(\lambda)-I)=0$. The  lemma below permits us to find $\nu$ as
the maximal $\lambda$ being the solution to equation
$\det\widetilde{D}(\lambda)=0$ or the equivalent equation
$\det\widehat{D}(\lambda)=0$. Observe that the entries of matrix
$D(\lambda)$ are expressed via the Laplace transforms of hitting
times \emph{under taboo sets} $W_i$, $i=1,\ldots,N$, whereas the
elements of $\widetilde{D}(\lambda)$ are represented in terms of the
Laplace transforms of hitting times \emph{without taboo}. At last,
the entries of $\widehat{D}(\lambda)$ involve the Laplace transforms
of transition probabilities only.
\begin{Lm}\label{L:Green_functions}
For any $\lambda>0$ one has $\det(D(\lambda)-I)=0$ if and only if
$\det\widetilde{D}(\lambda)=0$ or, equivalently, when
$\det\widehat{D}(\lambda)=0$. Moreover, these relations are true
even for $\lambda=0$ whenever $\eta$ is transient.
\end{Lm}

{\sc Proof.} Introduce matrices
$\widetilde{R}(\lambda)\!=\!\left(\widetilde{r}_{i,j}(\lambda)\right)_{i,j=1}^N$
and
$\widehat{R}(\lambda)\!=\!\left(\widehat{r}_{i,j}(\lambda)\right)_{i,j=1}^N$,
$\lambda>0$, with the corresponding entries
$\widetilde{r}_{i,j}(\lambda):=\delta_{i,j}+\left(1-\delta_{i,j}\right)F^\ast_{w_i,w_j}(\lambda)$
and ${\widehat{r}_{i,j}(\lambda):=G_{\lambda}(w_i,w_j)}$. By formula
(\ref{from_Chung}), the $j$-th column, $j=1,\ldots,N$, of matrix
$\widehat{R}(\lambda)$ is obtained from the $j$-th column of
$\widetilde{R}(\lambda)$ by multiplying it with
$G_{\lambda}(w_j,w_j)$ where $G_{\lambda}(w_i,w_j)>0$ in view of
irreducibility of $\eta$. Hence, $\det\widetilde{R}(\lambda)=0$ iff
$\det\widehat{R}(\lambda)=0$, for any $\lambda>0$. Let us show that
$\det\widehat{R}(\lambda)\neq0$ whenever $\lambda>0$. Similarly to
the arguments of Theorem~8 in \cite{Chung_60}, Ch.2, Sec.11, one can
derive that
\begin{equation}\label{hitting_time_new_t}
F_{w_i,w_j}(t)={_{W_j}F_{w_i,w_j}(t)} +\sum_{k=1,\,k\neq
j}^N\int\nolimits_0^t{_{W_k}F_{w_i,w_k}(t-u)}\,dF_{w_k,w_j}(u),
\end{equation}
for each $i,j=1,\ldots,N$. Note that the latter relation has a
natural interpretation. Namely, the path of the Markov chain from
$w_i$ to $w_j$ can either avoid the set $W_j$ or can hit some site
from $W_j$ and afterwards reach $w_j$. Applying the
Laplace-Stieltjes transform to (\ref{hitting_time_new_t}) we get
\begin{equation}\label{hitting_time_new}
F^\ast_{w_i,w_j}(\lambda)={_{W_j}F^\ast_{w_i,w_j}(\lambda)}
+\sum_{k=1,\,k\neq
j}^N{_{W_k}F^\ast_{w_i,w_k}(\lambda)}F^\ast_{w_k,w_j}(\lambda).
\end{equation}
Multiplying each side of equality (\ref{hitting_time_new}) by
$G_{\lambda}(w_j,w_j)$ and taking into account (\ref{from_Chung}) we
see that
$$G_{\lambda}(w_i,w_j)-\frac{\delta_{i,j}}{\lambda-q(w_i,w_i)}=
\sum_{k=1}^N{_{W_k}F^\ast_{w_i,w_k}(\lambda)}G_{\lambda}(w_k,w_j),
\quad i,j=1,\ldots,N.$$ Rewrite these identities in the following
matrix form
$$\widehat{R}(\lambda)-T(\lambda)=U(\lambda)\widehat{R}(\lambda)\quad\mbox{or, equivalently,}\quad
\widehat{R}(\lambda)\left(I-U(\lambda)\right)=T(\lambda)$$ where
$$T(\lambda)\!:=\!\left(\delta_{i,j}\left(\lambda-q(w_i,w_i)\right)^{-1}\right)_{i,j=1}^N,
\;\;\;
U(\lambda)\!:=\!\left({_{W_j}F^\ast_{w_i,w_j}(\lambda)}\right)_{i,j=1}^N,\;\;\lambda>0.
$$
Whence we deduce that $\det\widehat{R}(\lambda)\neq0$
for each $\lambda>0$, since $\det T(\lambda)$ is strictly positive.

It follows from (\ref{hitting_time_new}) that
$\left(D(\lambda)-I\right)\widetilde{R}(\lambda)=\widetilde{D}(\lambda)$.
Considering the determinants of the matrices at the left-hand and
the right-hand sides of the latter equality we come to the first
assertion of Lemma~\ref{L:Green_functions}. Its second assertion is
implied by representation
$\widetilde{d}_{i,j}(\lambda)G_{\lambda}(w_j,w_j)=\widehat{d}_{i,j}(\lambda)$
valid for each $\lambda>0$ and $i,j=1,\ldots,N$ by virtue of
(\ref{from_Chung}).

All the above reasoning holds true even for $\lambda=0$ whenever
$\eta$ is transient. Thus, Lemma~\ref{L:Green_functions} is proved
completely. $\square$

Let us consider some applications of our results to the models
studied earlier by different researchers.

\emph{Example 1.} Focus on a catalytic branching random walk on
${\bf Z}^d$, $d\in{\bf N}$, proposed in \cite{VTY}. This model is a
particular case of CBP if we set $\eta$ to be a symmetric and
space-homogeneous random walk on ${\bf Z}^d$ with a finite variance
of jump sizes as well as put $N=1$, $w_1={\bf 0}\in{\bf Z}^d$,
$\alpha_1=\alpha$, $\beta_1=1$ and ${f_1(\cdot)=f(\cdot)}$. Then we
deduce the same criticality condition as used in \cite{TV_SAM_13},
i.e. $D=1$ or, equivalently, ${\alpha f'(1)+(1-\alpha)\left(1+q({\bf
0},{\bf 0})^{-1}G^{-1}_0({\bf 0},{\bf 0})\right)=1}$. For recurrent
Markov chain $\eta$, $G^{-1}_0({\bf 0},{\bf 0})$ is assumed to be
$0$. Applying Theorem~\ref{T:MACBP} to catalytic branching random
walk on ${\bf Z}^d$ we come  to Theorem~4.1 and Theorem~4.2 in \cite{Y_TPA_11}
as well as some statements from Theorem~5 in \cite{TV_SAM_13} and
Theorem~1 in \cite{B_JTP_13}.

\emph{Example 2.} Consider catalytic branching process with a single
catalyst (located, say, at some site $w\in S$) studied in
\cite{DR_13}. Here the underlying motion of particles is governed by
an irreducible Markov chain $\eta$. Thus, such setting is less
restrictive than that in \emph{Example 1}. As shown in \cite{DR_13},
the asymptotic behavior of total and local particles numbers is
determined by the mean offspring number, produced by a particle at
the presence of the catalyst, being less than, equal to or greater
than $G^{-1}_0(w,w)+1$, and transience/recurrence of $\eta$. We come
to the same classification letting $N=1$, $w_1=w$,
$\alpha_1=\left(1-q(w,w)\right)^{-1}$, $\beta_1=1-q(w,w)$ and
$f'_1(1)=m$ in our CBP. Then, in view of (\ref{from_Chung}) and
evident  formula
$F_{w,w}(\infty)=\lim\nolimits_{\lambda\to0+}F^\ast_{w,w}(\lambda)$,
the value $sign(D-1)$ coincides with $sign(m-G^{-1}_0(w,w)-1)$.
Stress that in contrast to \cite{DR_13} we do not assume the
existence of all moments of the offspring number, i.e. finiteness of
$f^{(k)}_1(1)$ for any $k\in{\bf N}$. Some assertions of
Theorem~\ref{T:MACBP} in the case of CBP with a single catalyst are
stronger than the corresponding statements of Theorem~1 in
\cite{DR_13} (for instance, cf. relations (\ref{T:MACBP,sub_small})
and (\ref{T:MACBP,sub_big}) in our Theorem~\ref{T:MACBP} and point
iii)a) in Theorem~1 in \cite{DR_13}). However, some statements of
Theorem~1 in \cite{DR_13} are not covered by our
Theorem~\ref{T:MACBP} because they involve asymptotic estimates for
the moments of particles numbers in terms of local times at level
$w$ of the Markov chain $\eta$. The authors of \cite{DR_13} do not
discuss the asymptotic behavior of those local times.

\emph{Example 3.} Concentrate on the branching random walk on ${\bf
Z}^d$, $d\in{\bf N}$, with several sources investigated in
\cite{Y_MPS_13}. This model is a particular case of CBP such that
$\eta$ is a symmetric and space-homogeneous random walk on ${\bf
Z}^d$ with a finite variance of jump sizes and the symmetry of the
random walk fails only at a finite set of points of ${\bf Z}^d$. As
established in \cite{Y_MPS_13}, the rate of exponential growth of
the particles numbers is the maximal \emph{positive} $\lambda$
satisfying equation $\det\widehat{D}(\lambda)=0$ (this agrees with
our Lemma~\ref{L:Green_functions}). However, the necessary and
sufficient conditions of existence of this \emph{positive} solution
$\lambda$ were not found as noted in Conclusion in \cite{Y_MPS_13}.
It is worthwhile to remark that such necessary and sufficient
condition is provided in our present paper and by virtue of
Lemma~\ref{L:Green_functions} and Theorem~\ref{T:MACBP} it is just
$\rho(D)>1$, i.e. in our terms amounts to a supercritical regime of CBP.

Concluding the paper, we would like to observe that our approach of
combination of hitting times under taboo and auxiliary multi-type
Bellman-Harris processes permits us to obtain and justify the
effective classification of catalytic branching processes with
multiple catalysts. The  results of this paper are valid under
minimal restrictions on the character of motion and breeding of
particles. Thus they generalize
the previous works on closely related subject.  The developed approach can be
employed in the subsequent study of
other characteristics of catalytic branching process.

\section{Acknowledgements}\label{s:THANKS}$\mbox{ }$

The author is grateful to Professor M.A.Lifshits and Professor
V.A.Vatutin for useful discussions.

\end{document}